\documentclass{article}

\usepackage{fancyhdr}
\usepackage[T1]{fontenc}
\usepackage[english]{babel}
\usepackage[all]{xy}
\usepackage{amsfonts}
\usepackage{amsmath}
\usepackage{amssymb}
\usepackage{latexsym}
\usepackage{mathrsfs}  
\usepackage{graphicx} 
\usepackage[toc,page]{appendix}

\newtheorem{thm}{\bf Theorem}[section]
\newtheorem{cor}[thm]{\bf Corollary}

\newtheorem{lem}[thm]{\bf Lemma}
\newtheorem{prop}[thm]{\bf Proposition}
\newtheorem{defn}[thm]{\bf Definition}

\newtheorem{rem}[thm]{\bf Remark}
\newtheorem{rems}[thm]{\bf Remarks}
\newtheorem{exmp}[thm]{\bf Example}
\newtheorem{exmps}[thm]{\bf Examples}

\newcommand{\field}[1]{\mathbb{#1}}
\newcommand{\C}{\field{C}}
\newcommand{\R}{\field{R}}
\newcommand{\Q }{\field{Q}}
\newcommand{\Z }{\field{Z}}

\def\proof{{\parindent0pt {\bf Proof.\ }}}

\def\id{{\rm id}}

\def\F{\mathcal{F}_{\rm C}}
\def\P{\mathcal{P}_{\rm C}}
\def\I{\mathcal{I}_{\rm C^+}}
\def\C{\mathcal{C}_{\rm C}}
\def\SC{\mathcal{SC}_C}

\def\proj{\mathcal{P}}
\def\flat{\mathcal{F}}
\def\inj{\mathcal{I}}
\def\cot{\mathcal{C}}

\def\Mod{{\rm Mod\ }}

\def\Y{\mathcal{Y}}
\def\X{\mathcal{X}}

\def\A{\mathcal{A}}
\def\B{\mathcal{B}}
\def\Q{\mathcal{Q}}
\def\R{\mathcal{R}}
\def\W{\mathcal{W}}

\def\Gflat{{\rm \mathcal{GF}}}
\def\Gproj{{\rm \mathcal{GP}}}
\def\Ginj{{\rm \mathcal{GI}}}

\def\Gcot{{\rm \mathcal{GC}}}

\def\res{{\rm res}}
\def\cores{{\rm cores}}

\def\GC{{\rm G}_{\rm C}}
\def\GCF{{\rm G}_{\rm C}{\rm F}}
\def\GCC{{\rm G}_{\rm C}{\rm C}}

\def\Im{{\rm Im}}
\def\Coker{{\rm Coker}}
\def\Ker{{\rm Ker}}
\def\Add{{\rm Add}}
\def\add{{\rm add}}

\def\Ho{{\rm Ho}}

\def\Ext{{\rm Ext}}
\def\Tor{{\rm Tor}}
\def\Hom{{\rm Hom}}
\def\End{{\rm End}}
\def\Prod{{\rm Prod}}

\newcommand{\cqfd}
{\hspace{1cm}
\rule{2mm}{2mm}%
\medbreak%
\par%
}

\begin{document}

\title{Relative Gorenstein flat modules and Foxby classes and their model structures}

\author{Driss Bennis$^1$ \hskip 2cm  Rachid El Maaouy$^{2,}\footnote{Corresponding author}$ \\ \\ J. R. Garc\'{\i}a Rozas$^3$ \hskip 1,5cm Luis Oyonarte$^4$}

\date{}

\maketitle

\small{1: CeReMaR Research Center, Faculty of Sciences, B.P. 1014, Mohammed V University in Rabat, Rabat, Morocco.
	
	\noindent e-mail address: driss.bennis@um5.ac.ma; driss$\_$bennis@hotmail.com
	
	2: CeReMaR Research Center, Faculty of Sciences, B.P. 1014, Mohammed V University in Rabat, Rabat, Morocco.
	
	\noindent e-mail address: rachid\_elmaaouy@um5.ac.ma; elmaaouy.rachid@gmail.com
	
	3: Departamento de  Matem\'{a}ticas, Universidad de Almer\'{i}a, 04071 Almer\'{i}a, Spain.
	
	\noindent e-mail address: jrgrozas@ual.es
	
	4: Departamento de  Matem\'{a}ticas, Universidad de Almer\'{i}a, 04071 Almer\'{i}a, Spain.
	
	\noindent e-mail address: oyonarte@ual.es}

\begin{abstract}
A model structure on a category is a formal way of introducing a homotopy theory on that category, and if the model structure is abelian and hereditary, its homotopy category is known to be triangulated. So, a good way to both build and model a triangulated category is to build a hereditary abelian model structure. 

Given a ring $R$ and a (non-necessarily semidualizing) left $R$-module $C$, we introduce and study new concepts of relative Gorenstein cotorsion and cotorsion modules: $\GC$-cotorsion and (strongly) $\C$-cotorsion. As an application, we prove that there is a unique hereditary abelian model structure on the category of left $R$-modules, in which the cofibrations are the monomorphisms with $\GC$-flat cokernel and the fibrations are the epimorphisms with $\C$-cotorsion kernel belonging to the Bass class $\B_C(R)$. In the second part, when $C$ is a semidualizing $(R,S)$-bimodule, we investigate the existence of abelian model structures on the category of left (resp., right) $R$-modules where the cofibrations are the epimorphisms (resp., monomorphisms)  with kernel (resp., cokernel) belonging to the Bass (resp., Auslander) class $\B_C(R)$ (resp., $\A_C(R)$).

We also study the class of $\GC$-flat modules and the Bass class from the Auslander-Buchweitz approximation theory point of view. We show that they are part of weak AB-contexts. As the concept of weak AB-context can be dualized, we also give dual results that involve the class of $\GC$-cotorsion modules and the Auslander class. 

\end{abstract}



\section{Introduction}

Model categories were first introduced by Quillen in \cite{Qui67}. Roughly speaking, a model structure on a bi-complete category $\A$ is given by three classes of morphisms of $\A$, called cofibrations, fibrations and weak equivalences, from which it is possible to introduce a homotopy theory in $\A$. Model categories are interesting because they establish the theoretical framework for formally inverting the weak equivalences. We refer the reader to \cite{Hov99} for the definition and basic results. 

Later, Hovey (\cite{Hov02}), and also  Beligiannis and  Reiten (\cite[Chapter VIII]{BR07}), proposed the notion of an abelian model structure on an abelian category $\A$, as a model structure that is compatible with the abelian context of the category. Hovey proved (\cite[Theorem 2.2]{Hov02}) that there is a one-to-one correspondence between the class of abelian model structures and the class of Hovey triples (i.e., three classes of objects satisfying some properties). This establishes a relation between the theory of model categories and the representation theory via cotorsion pairs. If the abelian model structure is hereditary (i.e., the two complete cotorsion pairs induced by the Hovey triple are hereditary), then its homotopy category, $\Ho(\A)$, is the stable category of a Frobenius category, and so it is triangulated. This is an important situation in which one obtains a triangulated category from the point of view of model category theory. We refer the reader to \cite{Gil16} for a recent and interesting survey on hereditary abelian model structures. 

The first purpose of this paper is to construct new  hereditary abelian model structures on the category of left $R$-modules that involve the class of $\GC$-flat modules and the well-known Foxby classes, and use the homotopy categories of these model structures to further study these classes. A second purpose of this paper is to continue the study initiated in \cite{BDGO21,BEGO21b,BEGO22,BGO16} of relative Gorenstein homological algebra with respect to a non-necessary semidualizing module $C$. In particular, we will prove results that are new even in the case where $C$ is a semidualizing module.

The study of the Gorenstein homological dimension with respect to a semidualizing module $C$ goes back to Golod \cite{Gol84} when he introduced the finitely generated modules of finite $\GC$-dimension over a commutative noetherian ring $R$, as a relative theory to that of Auslander and Bridger \cite{AB69}. Holm and  J\symbol{'370}rgensen in \cite{HJ06} extended Golod's study of $\GC$-dimension of finitely generated modules to the case of arbitrary modules over commutative noetherian rings. They considered the $C$-Gorenstein projective  and $C$-Gorenstein flat $R$-modules. These modules are the analogues of modules of $\GC$-dimension $0$. Recently, these approaches have been unified in \cite{BGO16,BEGO22} with the use of $\GC$-projective $R$-modules and $\GC$-flat $R$-modules. 

We will also be interested in the case where the class of $\GC$-flat modules and the Bass class are part of weak Auslander-Buchweitz contexts in Hashimoto's terminology \cite{Has20}. This allows us, among other things, to deduce the existence of certain Auslander-Buchweitz approximations for $R$-modules of finite $\GC$-flat and $\B_C$-injective dimensions, respectively.

This article is organized as follows. 

In Section 2 we give the key definitions and preliminary results necessary for the rest of the paper. 

In Section 3 we introduce and study two new concepts related to relative cotorsion modules: strongly $\C$-cotorsion and $n$-$\C$-cotorsion modules for a given integer $n\geq 1$. We are mainly interested in their links with other known classes of modules such as cotorsion modules (Proposition \ref{Coto and C-cot} and Corollary \ref{Coto and strong C-cot}), as well as in their homological properties. It is investigated when these new classes are the right side class of a (perfect, complete, hereditary) cotorsion pair (Theorem \ref{cot pair-relative flat}).

In Section 4 we introduce and investigate a new concept of Gorenstein cotorsion modules: $\GC$-cotorsion.  We characterize when the pair ($\GC$-flat,$\GC$-cotorsion) is a hereditary and perfect cotorsion pair (Theorem \ref{(GwF,GwC) cot pair}). Then, two different descriptions of the core of this cotorsion pair is given (Proposition \ref{(Fc,Cc) core} and Proposition \ref{2nd descr of the core of (GcF,GcC)}). As applications, we show that $\GCF(R)$, the class of $\GC$-flat $R$-modules, is part of a left weak AB-context and that $\GCC(R)$, the class of $\GC$-cotorsion $R$-modules, is part of a right weak AB-context (Theorem \ref{Gc-flat weak AB-cont}). As a second application, we construct a new hereditary abelian model structure on the category of left $R$-modules, where the cofibrations are precisely the monomorphisms with $\GC$-flat cokernel and the fibrations are the epimorphisms with $\C$-cotorsion kernel belonging to $\B_C(R)$, the Bass class (Theorem \ref{GwF model structure}).

In Section 5 we investigate when $\B_C(R)$ is part of a right weak AB-context. Next, we find necessary and sufficient conditions to guarantee the existence of a hereditary abelian model structure whose fibrations are epimorphisms with kernel lying in $\B_C(R)$ and whose cofibrations are  monomorphisms with cokernel lying in $^\perp\Add_R(C)$ (Theorem \ref{bass-mod-struc}). At the end of this section we relate these model structures to the Frobenius category relative to $C$, the so called $C$-Frobenius category, recently introduced in \cite{BEGO21a}. It is shown that when the category of left $R$-modules is $C$-Frobenius, the homotopy category is an extension of the stable model category of a Frobenius ring. Of course, dual results concerning the Auslander class are investigated.

\section{Preliminaries}

Throughout this paper, $R$ and $S$ will be  associative (non-necessarily commutative) rings with identity, and all modules will be, unless otherwise specified, unital left $R$-modules or left $S$-modules. When right $R$-modules need to be used, they will be denoted as $M_R$, while in these cases left $R$-modules will be denoted by $_R M$. 
We use $\inj(R)$, $\proj(R)$, $\flat(R)$, and $\cot(R)$ to denote the classes of injective, projective, flat, and cotorsion $R$-modules, respectively.

From now on $C$ will stand for an $R$-module, $S$ for its endomorphism ring, $S=\End_R(C)$, and $C^+$  for the character module of $C$, $C^+=\Hom_{\Z}(C,\mathbb{Q}/\Z)$. Notice that $C$ is an $(R,S)$-bimodule and $C^+$ an $(S,R)$-bimodule. We use $\Add_R(C)$ (resp., $\Prod_R(C)$) to denote the class of all $R$-modules which are isomorphic to direct summands of direct sums (resp., direct products) of copies
of $C$, and we write $\add_R(C)$ when such a direct sum is finite.

Given a class  $\X$ of $R$-modules and a class $\mathcal{Y}$ of right $R$-modules, an $\mathcal{X}$-resolution of an $R$-module $M$ is an exact complex $$\cdots\to X_1\to X_0\to M\to 0$$ where $X_i\in\mathcal{X}$. An $\mathcal{X}$-coresolution of $M$ is defined dually. A complex $\textbf{X}$ in $R$-Mod is called $\left(\mathcal{Y}\otimes_R-\right)$-exact (resp., $\Hom_R(\mathcal{X},-)$-exact, $\Hom_R(-,\mathcal{X})$-exact) if $Y\otimes_R\textbf{X}$  (resp., $\Hom_R(X,\textbf{X})$, $\Hom_R(\textbf{X},X)$) is an exact complex for every $Y\in\mathcal{Y}$ (resp., $X\in\mathcal{X}$).

A class $\mathcal{X}$ of $R$-modules is projectively resolving if the class $\mathcal{X}$ is closed under extensions, kernels of epimorphisms, and $\proj(R)\subseteq\mathcal{X}$. The class $\X$  is left thick if the class $\mathcal{X}$ is closed under extensions, kernels of epimorphisms, and direct summands. Injectively coresolving and right thick classes are defined dually. $\mathcal{X}$ is thick if it is left and right thick.

\medskip
\noindent\textbf{Cotorsion pairs.}
Given a class of $R$-modules $\X$  and an integer $n\geq 1$, we use the following standard notations:  $$\mathcal{X}^{\perp_n}=\{M\in R\text{-Mod}|\Ext^i_R(X,M)=0, \forall X\in\mathcal{X},\forall i=1,\dots ,n \}.$$ $$^{\perp_n}\mathcal{X}=\{M\in R\text{-Mod}|\Ext^i_R(M,X)=0, \forall X\in\mathcal{X},\forall i=1,\dots ,n \}.$$ In particular, we set $\mathcal{X}^{\perp}=\mathcal{X}^{\perp_1}$, $^{\perp}\mathcal{X}=$ $^{\perp_1}\mathcal{X}$, $\mathcal{X}^{\perp_\infty}=\cap_{n\geq 1}\mathcal{X}^{\perp_n}$ and   $^{\perp_\infty}\mathcal{X}=\cap_{n\geq 1}\;^{\perp_n}\mathcal{X}$.

Given a class of $R$-modules $\mathcal{X}$, an $\mathcal{X}$-precover of a module $M$ is a morphism $f:X\to M$ with $X\in \mathcal{X}$, in such a way that $f_*:\Hom_R(X',X)\to \Hom_R(X',M)$ is surjective for every $X'\in \mathcal{X}$. The $\mathcal{X}$-precover $f$ is said to be special if it is an epimorphism and $\Ker (f) \in \X^{\perp}$, and it is said to be an $\mathcal{X}$-cover provided that every endomorphism $g:X\to X$ such that $fg=f$ is an automorphism of $X$. If every module has a $\mathcal{X}$-precover, a special $\X$-precover or a $\X$-cover, then the class $\mathcal{X}$ is said to be precovering, special precovering or covering, respectively. Special $\mathcal{X}$-preenvelopes, preenvelopes and envelopes can be defined dually.

A pair $\left(\mathcal{A},\mathcal{B}\right)$ of classes of modules is called a cotorsion pair if $\mathcal{A}^\perp=\mathcal{B}$ and $\mathcal{A}=$ $^\perp\mathcal{B}$. If $\left(\mathcal{A},\mathcal{B}\right)$ is a cotorsion pair, then the class $\mathcal{A}\cap\mathcal{B}$ is called the core of the cotorsion pair. A cotorsion pair $\left(\mathcal{A},\mathcal{B}\right)$ is said to be hereditary if $\mathcal{A}$ is projectively resolving, or equivalently, if $\mathcal{B}$ is injectively coresolving.  A cotorsion pair $\left(\mathcal{A},\mathcal{B}\right)$ is said to be complete if any $R$-module has a special $\B$-preenvelope, or equivalently, any $R$-module has a special $\A$-precover. A cotorsion pair $\left(\mathcal{A},\mathcal{B}\right)$ is said to be perfect if every module has an $\mathcal{A}$-cover and a $\mathcal{B}$-envelope. It is clear by its definition that a perfect cotorsion pair $(\mathcal{A},\mathcal{B})$ is complete. The converse holds when $\mathcal{A}$ is closed under direct limits \cite[Corollary 2.3.7]{GT12}. A cotorsion pair $\left(\mathcal{A},\mathcal{B}\right)$ is said to be cogenerated by a set if there is a set $\mathcal{X}$ such that $\mathcal{B}=\mathcal{X}^\perp$. It is a well known fact that any cotorsion pair cogenerated by a set is complete (see \cite[Theorem 7.4.1]{EJ00} for example).

The following two lemmas will be used in several places in the paper. The second lemma is known as Wakamatsu's Lemma (see \cite[Lemmas 2.1.1 and 2.1.2]{Xu96}).

\begin{lem}\label{gen-cog sets}
Let $\mathcal{X}$ be a set of $R$-modules.
\begin{enumerate}
\item $\mathcal{X}^{\perp_\infty}=M^\perp$ for some $R$-module $M$.
\item $^{\perp_\infty}\mathcal{X}=$ $^\perp M$ for some $R$-module $M$.
\end{enumerate}
\end{lem}
\proof The proof of (2) is dual to that of (1), so we only prove (1). 

Let $X$ be the direct sum of all the modules in $\mathcal{X}$, consider any projective resolution of $X$, $$\cdots\to P_1\stackrel{f_1}{\to} P_0\stackrel{f_0}{\to} X\to 0,$$ and let $K_{i+1}=\Ker(f_i)$. Clearly, for any $i\geq 1$ and any $R$-module $A$, we have  $\Ext^1_R(K_i,A)\cong \Ext^{i+1}_R(X,A)$. If we let $M=X\oplus (\oplus_{i\geq 1}K_i)$, then 
$\Ext^1_R(M,A)\cong \Ext^1_R(X,A)\oplus (\prod_{i\geq 1}\Ext^1_R(K_i,A))\cong\prod_{i\geq 1}\Ext^{i}_R(X,A). \text{\cqfd}$

\begin{lem} [\cite{Xu96}, Lemmas 2.1.1 and 2.1.2] \label{Wakamtsu Lemma}
Assume that  $\mathcal{X}\subseteq R$-\Mod is closed under extensions. 
\begin{enumerate}
\item If $\varphi :X\to M$ is an epic $\mathcal{X}$-cover of $M$, then $\varphi$ is a special $\mathcal{X}$-precover. 
\item  If $\varphi :M\to X$ is a monic $\mathcal{X}$-envelope of $M$, then $\varphi$ is a special $\mathcal{X}$-preenvelope. 
\end{enumerate}
\end{lem}

Consequently, every covering class $\mathcal{X}\subseteq R$-\Mod closed under extensions and such that $\proj(R)\subseteq \mathcal{X}$ is special precovering. Dually, every enveloping class $\mathcal{X}\subseteq R$-\Mod closed under extensions and such that $\inj(R)\subseteq \mathcal{X}$ is special preenveloping.

\medskip\noindent\textbf{Abelian model structures.}
Given a bicomplete category $\A$, a model structure on $\A$ is given by three classes of morphisms of $\A$ called cofibrations, fibrations and weak equivalences, satisfying a set of axioms (see \cite[Definition 1.1.3.]{Hov99}). If $\mathcal{M}$ is a model structure on $\A$, then it  provides a general framework to study homotopy theory. By this we mean that we can associate to $\mathcal{M}$ the homotopy category, $\Ho_\A(\mathcal{M})$, which  is defined by formally inverting the weak equivalences of $\mathcal{M}$. More precisely, $\Ho_\A(\mathcal{M})$ is obtained after localizing $\A$ at the class of weak equivalences (see \cite[Sec. 1.2]{Hov99} for more details).

Hovey restricted the study of model structures to abelian categories. He defined in \cite{Hov02} the notion of an abelian model structure and showed in \cite[Theorem 2.2]{Hov02} that there is a close link between some model structures on an abelian category and the cotorsion pairs in this category. Namely, an abelian model structure on any abelian category $\A$, is equivalent to a triple $\mathcal{M}=(\mathcal{Q},\mathcal{W},\mathcal{R})$ of classes of objects in $\cot$ such that $\mathcal{W}$ is thick and $(\mathcal{Q},\mathcal{W}\cap\mathcal{R})$ and $(\mathcal{Q}\cap\mathcal{W},\mathcal{R})$ are complete cotorsion pairs. In this case, $\Q$ is precisely the class of cofibrant objects, $\R$ is precisely the class of fibrant objects, and $\W$ is the class of trivial objects of the model structure. Here, an abelian model structure in the sense of \cite[Definition 1.1.3]{Hov99} is one in which: (a) a morphism is a (trivial) cofibration if and only if it is a monomorphism with (trivially) cofibrant cokernel; and (b) a morphism is a (trivial) fibration if and only if it is an epimorphism with (trivially) fibrant kernel.

Hovey's correspondence shows that  we can shift all our attention from morphisms (cofibrations, weak equivalences, and fibrations) to objects (cofibrant, trivial, and fibrant). In this paper, we identify any abelian model structure with such a triple and we often call it a Hovey triple.

A Hovey triple is hereditary if the two associated cotorsion pairs are hereditary. In general, it is difficult to prove that a particular category has a model structure. However, Gillespie \cite{Gil15} gave a new and useful method for constructing hereditary abelian model structures as follows.

\begin{thm}[\cite{Gil15}, Theorem 1.1]\label{Hovey's correpondence}
Given two complete and hereditary cotorsion pairs $(\widetilde{\Q},\mathcal{R})$ and $(\mathcal{Q},\widetilde{\R})$ in an abelian category $\mathcal{A}$ such that  $\widetilde{\R}\subseteq \R$ (or equivalently, $\widetilde{\Q}\subseteq \Q$) and  $\widetilde{\Q}\cap \mathcal{R}=\mathcal{Q}\cap\widetilde{\mathcal{R}}$, there is a unique thick class $\mathcal{W}$ such that $(\mathcal{Q},\mathcal{W},\mathcal{R})$ is a Hovey triple. Moreover, the class $\mathcal{W}$ is characterized by:
\begin{eqnarray*}
\W=\{X\in\A|\exists\text{ an exact  sequence } 0\to X\to R'\to Q'\to 0\text{ with } R'\in\mathcal{\widetilde{R}}, Q'\in\mathcal{\widetilde{Q}}\}\\
=\{X\in\A|\exists\text{ an exact  sequence } 0\to R'\to Q'\to X\to 0\text{ with } R'\in\mathcal{\widetilde{R}}, Q'\in\mathcal{\widetilde{Q}}\}.\end{eqnarray*}
\end{thm}

\medskip\noindent\textbf{Weak AB-contexts.} Given a full subcategory $\X$ of an abelian category $\A$, we denote by $\res(\widehat{\X})$ the full subcategory of objects in $\A$ having a finite $\X$-resolution. Dually, we denote by $\cores(\widehat{\X})$ the full subcategory of objects in $\A$ having a finite $\X$-coresolution. 

Following \cite{BMPS19}, a triple $(\X,\Y,\omega)$ of full subcategories of $\A$ is said to be a left weak Auslander-Buchweitz context (or a left weak AB-context for short) if the following three conditions are satisfied:

(AB1) $\mathcal{X}$ is left thick.

(AB2) $\mathcal{Y}$ is right thick and $\mathcal{Y}\subseteq \res(\widehat{\mathcal{X}})$.

(AB3) $\omega=\X\cap\Y$ and $\omega$ is an injective cogenerator for $\X$, that is, $\omega \subseteq \X\cap \X^{\perp_\infty}$ and for any $X\in\X$ there exists an exact sequence $0\to X\to W\to X'\to 0$ such that $X'\in\X$ and $W\in\omega$.

Right weak AB-contexts can be defined dually, and these were considered in \cite{GMZ18} under the name weak co-AB-contexts.

Becerril, Mendoza, P\'{e}rez and Santiago introduced in \cite{BMPS19} the notion of left Frobenius pairs and showed that there is a one-to-one correspondence between the class
of left weak AB-contexts and that of Frobenius pairs. Recall that a pair $(\X,w)$ of full subcategories of $\A$ is a left Frobenius pair if $\X$ is left thick, and $w$ is closed under direct summands and an injective cogenerator for $\mathcal{X}$. Right Frobenius pairs $(\omega,\Y)$ can be defined dually.

So, looking for left weak AB-contexts is the same as looking for left Frobenius pairs. The following two lemmas are the key to find our left and right weak AB-contexts in this  paper. 

\begin{lem}\label{left Frobnius}
The following assertions hold:
\begin{enumerate}
\item (\cite[Theorem 5.4(1)]{BMPS19}) If $(\mathcal{X},\omega)$ is a left Frobenius pair in $\A$, then $(\mathcal{X},\res(\widehat{\omega}), \omega)$ is a left weak AB-context. Conversely, if $(\mathcal{X},\mathcal{Y},\omega)$ is a left weak AB-context in $\A$, then $(\mathcal{X},\omega)$ is a left Frobenius pair.
\item  (\cite[Proposition 2.5]{LY22}) If $(\X,\Y)$ is a complete hereditary cotorsion pair in $\A$, then $(\X,\X\cap \Y)$ is a left Frobenius pair.
\end{enumerate}
\end{lem}

In this paper we will also be interested in right weak AB-contexts. A dual result of  Lemma \ref{left Frobnius} is therefore needed. For the reader's convenience we state it here. It is not mentioned in either \cite{BMPS19} or \cite{LY22}. However, its proof is simply a dual to that of \cite[Theorem 5.4(1)]{BMPS19} and \cite[Proposition 2.5]{LY22}.

\begin{lem}\label{right Frobnius}
The following statements hold:
\begin{enumerate}
\item If $(\omega,\Y)$ is a right Frobenius pair in $\A$, then $(\omega,\cores (\widehat{\omega}),\Y)$ is a right weak AB-context. Conversely, if $(\omega,\X,\Y)$ is a right weak AB-context in $\A$, then $(\omega,\Y)$ is a right Frobenius pair.
\item If $(\X,\Y)$ is a complete hereditary cotorsion pair in $\A$, then $(\X\cap\Y, \Y)$ is a right Frobenius pair. 
\end{enumerate}
\end{lem}

\medskip\noindent \textbf{Relative Gorenstein flat modules.} We recall the concept of $\GC$-flat $R$-modules and all related notions.

\begin{defn}[\cite{BEGO22}] Let $C$ be an $R$-module. An  $R$-module $M$ is said to be $\F$-flat (resp., $\mathcal{I}_C$-injective) if $M^+\in \Prod_R(C^+)$ (resp., $M\in \Prod_R(C)$). 

We denote the class of all $\F$-flat modules as $\F(R)$ and that of all $\mathcal{I}_C$-injective modules as $\mathcal{I}_C(R)$.
\end{defn}

\begin{rem}Note that an $R$-module $M$ is $\F$-flat if and only if $M^+$ is an $\I$-injective right $R$-module, and that the class of $\F$-flat modules is clearly closed under direct sums and direct summands.
\end{rem}

\begin{exmps}\label{exmps of Fc-flat} Here are some examples of $\F$-flat and $\I$-injective modules:
\begin{enumerate}
\item If $C$ is a flat generator of $R$-{\rm Mod}, then $$\F(R)=\flat(R) \text{ and } \Prod_R(C^+)=\inj(R).$$
\item (\cite[Proposition 3.3]{BEGO22}) If $_RC$ is finitely presented then $$ \F(R)=C\otimes_S\flat(S)  \text{ and } \Prod_R(C^+)=\Hom_S(C,\inj(S)).$$
\end{enumerate}
\end{exmps}

\begin{defn}[\cite{HW06}, Definition 2.1] An $(R,S)$-bimodule $C$ is semidualizing if the following conditions hold:
\begin{enumerate}
\item $_RC$ and $C_S$ both admit a degreewise finite projective resolution.
\item $\Ext_R^{\geq 1}(C,C)=\Ext_S^{\geq 1}(C,C)=0.$
\item The natural homothety maps $R\rightarrow\Hom_S(C,C)$ and $S \rightarrow \Hom_R(C,C)$ both are ring isomorphisms.
\end{enumerate}
\end{defn}

\begin{rem} 
Without any assumption on $C$ we always have that every $C$-flat module (resp., $C$-injective right $R$-module), in the sense of \cite[Definition 5.1]{HW06}, is $\F$-flat (resp., $\I$-injective). However, when $C$ is semidualizing the classes of $\F$-flat and $\I$-injective modules coincide with those of $C$-flat and $C$-injective modules, respectively.
\end{rem}

\begin{defn}[\cite{BGO16}] An $R$-module $C$ is said to be w-tilting if it satisfies the following two properties:
\begin{enumerate}
\item $C$ is $\Sigma$-self-orthogonal, that is, $\Ext_R^{i\geq 1}(C,C^{(I)})=0$ for every set $I$.
\item There exists a $\Hom_R(-,\Add_R(C))$-exact $\Add_R(C)$-coresolution $$\textbf{X}:0\rightarrow R\rightarrow C_{-1}\rightarrow C_{-2}\rightarrow\cdots$$
\end{enumerate}
Dually, we define $\prod$-self-orthogonal and w-cotilting modules.
\end{defn}

The notion of a $w^+$-tilting module was recently introduced in \cite{BEGO22}, where the authors studied Gorenstein flat modules with respect to a non-necessary semidualizing bimodule. It properly generalizes  the concepts of semidualizing, Wakamatsu tilting and w-tilting modules (see \cite[Proposition 4.3]{BEGO22} and \cite[Example 2.2]{BGO16}). 

\begin{defn} [\cite{BEGO22}, Definition 4.1] An $R$-module $C$ is said to be $w^+$-tilting if it satisfies the following two properties:
\begin{enumerate}
\item $C$ is $\prod$-$\Tor$-orthogonal, that is, $\Tor^R_{i\geq 1}((C^+)^I,C)=0$ for every set $I$.
\item There exists a  $\left( \Prod_R(C^+)\otimes-\right) $-exact $\F(R)$-coresolution $$\textbf{X}:0\rightarrow R\rightarrow C_{-1}\rightarrow C_{-2}\rightarrow\cdots$$
\end{enumerate}
\end{defn}

\begin{defn} [\cite{BEGO22}, Definition 4.4]\label{GcF-def} Let $C$ be an $R$-module. An $R$-module $M$ is said to be $\GC$-flat if there exists an exact and $\left(\Prod_R(C^+)\otimes_R-\right)$-exact sequence $$\textbf{X}:\cdots \rightarrow F_1\rightarrow F_0 \rightarrow C_{-1}\rightarrow C_{-2}\rightarrow\cdots$$ with each $C_i\in\F(R)$ and $F_j\in\flat(R)$, such that $M\cong \Im(F_0\to C_{-1})$. The $\left(\Prod_R(C^+)\otimes_R-\right)$-exact exact sequence $\textbf{X}$ is called a complete $\F$-flat resolution of $M$.
	
The class of all $\GC$-flat $R$-modules is denoted by $\GCF(R)$.
\end{defn}

Recall \cite{BGO16} that an $R$-module $M$ is said to be $\GC$-projective if there exists a $\Hom_R(-,\Add_R(C))$-exact exact sequence $\cdots \to P_1\to P_0\to C_{-1}\to C_{-2}\to \cdots$, where $P_i\in\proj(R)$, $C_j\in\Add_R(C)$, and such that $M\cong \Im(P_0\to C_{-1})$.

It is still an open question whether or not (or at least when) every $\GC$-projective module is $\GC$-flat (even when $C=R$). The following result will be used later. It gives a partial answer to this question.

\begin{prop} \label{Gc-proj is Gc-flat} Assume that $_RC$ has a degreewise finite projective resolution. The following statements hold:
\begin{enumerate}
\item $_RC$ is $\prod$-$\Tor$-orthogonal if and only if it is $\Sigma$-self-orthogonal.
\item  Assume that $_RC$ is $w$-tilting.  Then, every $\GC$-projective $R$-module which has a degreewise finite projective resolution is $\GC$-flat. In particular, $_RC$ is $w^+$-tilting. 
\end{enumerate} 
\end{prop}
\proof 1. Since $_RC$ has a degreewise finite projective resolution, this assertion follows from the isomorphism $\Tor^R_i((C^+)^{I},C)\cong \Ext_R^i(C,C^{(I)})^+$ for every set $I$ and $i\geq 1$, by \cite[Theorem 3.2.13 and Remark 3.2.27]{EJ00}.

2. Let $M$ be a $\GC$-projective $R$-modules admitting a degreewise finite projective resolution. A similar argument to that of \cite[Proposition 2.5]{BEGO21a} gives a $\Hom_R(-,C)$-exact exact sequence $$\textbf{X}: \cdots \to P_1\to P_0\to C_{-1}\to C_{-2}\to \cdots,$$ where the $P_i's$ are all finitely generated and projective, each $C_j\in\add_R(C)$, and $M\cong \Im(P_0\to C_{-1})$.

Now using \cite[Theorems 3.2.11 and 3.2.22]{EJ00}, we get isomorphisms of complexes $$(C^+)^{I}\otimes_R \textbf{X}\cong (C^+\otimes_R \textbf{X})^{I}\cong (\Hom_R(\textbf{X},C)^+)^{I}$$ for every set $I$. But $\Hom_R(\textbf{X},C)$ is exact so is $X\otimes_R\textbf{X}$ is exact for every $X\in\Prod_R(C^+)$. Hence, $M$ is $\GC$-flat.

Finally, to see that $C$ is $w^+$-tilting note that  $C$ being $w$-tilting means that $C$ and $R$ are $\GC$-projective. But since both $R$ and $C$ have a degreewise finite projective resolution, we get that they are also $\GC$-flat $R$-modules. Hence, $C$ is $w^+$-tilting. \cqfd

Certain properties of $\GC$-flat modules depend on the fact that the class $\GCF(R)$ is closed under extensions. But, unlike the class of $\GC$-projective and $\GC$-injective modules, it is unknown for what modules $C$ the class of $\GC$-flat modules is closed under extensions (see \cite[Question 4]{BEGO22}). If we assume that $_RC$ is $\prod$-$\Tor$-orthogonal, then there are some situations under which this closure property holds:

\begin{itemize}
\item[(a)] When $C$ is a flat generator of $R$-Mod. This situation follows from the equalities  $\F(R)=\flat(R)$ and $\GCF(R)=\Gflat(R)$ by \cite[Proposition 4.23]{BEGO22}  and the work of \v{S}aroch and  \v{S}\v{t}ov\'{\i}\u{c}ek  \cite[Theorem 4.11]{SS20}. 

\item[(b)] When every $R$-module has finite $\F$-flat dimension (see \cite[Corollary 7.5]{BEGO22}).

\item[(c)] When $\F(R)$ is closed under direct products (\cite[Corollary 4.13]{BEGO22}). In particular, when $S$ is right coherent and both $_RC$ and $C_S$ are finitely presented (see Lemma \ref{exmp of GwF-closed}).
\end{itemize}

This discussion naturally leads to the following question: 


\medskip
\textbf{Question:} Is the class of $\GC$-flat $R$-modules closed under extensions for any $\prod$-$\Tor$-orthogonal $R$-module $C$?

\medskip
Throughout, the following definition is needed.

\begin{defn}[\cite{BEGO22}, Definition 4.14] Let $C$ be an $R$-module. A ring $R$ is said to be $\GCF$-closed provided that the class of $\GC$-flat $R$-modules is closed under extensions.
\end{defn} 

\medskip\noindent \textbf{Foxby classes.} Associated to the bimodule $_RC_S$ we have the Auslander and Bass classes, $\A_C(S)$ and $\B_C(R)$, respectively, defined as follows:

$\bullet$ $\A_C(S)$ is the class of all $S$-modules $M$ satisfying:

$\Tor^S_{\geq 1}(C,M)=0=\Ext_R^{\geq 1}(C,C\otimes_S M)$ and the canonical map $$\mu_M:M\to \Hom_R(C,C\otimes_SM)$$ is an isomorphism of left $S$-modules.

$\bullet$ $\B_C(R)$ consists of all $R$-modules $N$ satisfying:

$\Ext_R^{\geq 1}(C,N)=0=\Tor^S_{\geq 1}(C,\Hom_R(C,N))=0$ and the canonical map $$\nu_N:C\otimes_S\Hom_R(C,N)\to N$$ is an isomorphism of $R$-modules.

On the other hand, one can define the classes $\A_C(R)$ and $\B_C(S)$ of right $R$-modules and right $S$-modules, respectively.

We will refer to the modules in $\A_C$ and $\B_C$ as $\A_C$-Auslander and $\B_C$-Bass modules, respectively. 

It is an important property of Auslander and Bass classes that they are equivalent under the following pair of functors (\cite[Proposition 2.1]{GLT15}): $$\xymatrixcolsep{5pc}\xymatrix{ \B_C(R)\ar@/^1pc/[r]^{\Hom_R(C,-)}  &\A_C(S) \ar@/^1pc/[l]^{C\otimes_S-} } \text{ and } \xymatrixcolsep{5pc}\xymatrix{ \A_C(R)\;\ar@/^1pc/[r]^{-\otimes_RC}  &\B_C(S) \ar@/^1pc/[l]^{\Hom_S(C,-)} }.$$

Consequently, Bass classes can be defined via Auslander  classes and vice-versa: $$\B_C(R)=C\otimes_S\A_C(S) \text{ and } \A_C(R)=\Hom_S(C,\B_C(S)).$$

Recall that an $R$-module $M$ is called self-small if the canonical morphism $$\Hom_R(M,M^{(I)})\to\Hom_R(M,M)^{(I)}$$ is an isomorphism, for every set $I$. Examples of self-small modules are finitely generated modules. Note that the module $_RC$ is self-small if and only if, for every set $I$, the canonical map $\mu_{S^{(I)}}:S^{(I)}\to \Hom_R(C,C\otimes_S  S^{(I)})$ is an isomorphism. 

Foxby classes are expected to satisfy certain properties under dual assumptions. Inspired by this duality, we propose the dual notion to that of self-small.

\begin{defn}An $R$-module $_RM$ is said to be self-co-small, if the canonical morphisms $$(M^+)^I\otimes_R M\to (M^+\otimes_R M)^I \text{ and } M^+\otimes_RM\to \Hom_R(M,M)^+$$ are isomorphisms for every set $I$. 	
\end{defn}

\begin{rem}
\item \begin{enumerate}
\item The module $_RC$ is self-co-small if and only if the canonical morphism $$\nu_{(S^+)^I}:\Hom_S(C,(S^+)^I)\otimes_RC\to (S^+)^I$$ is an isomorphism for every set $I$.
\item Any finitely presented $R$-module is self-co-small by \cite[Theorems 3.2.11 and 3.2.22]{EJ00}.
\end{enumerate}
\end{rem}
 
With respect to the terminology used in this paper, modules in the class $\Add_R(C)$ will be called $\P$-projective. By \cite[Proposition 3.1]{BDGO21}, when $C$ is self-small, the class of $\P$-projective modules coincides with that of $C$-projective modules, that is, modules in the class  $C\otimes_S \proj(S)$ (see \cite{HW06}).

It is straightforward to prove the following.

\begin{lem}\label{equiv of categ} When $_RC$ is self-small the pair of functors $$\xymatrixcolsep{5pc}\xymatrix{ \Add_R(C)\ar@/^1pc/[r]^{\Hom_R(C,-)} & \proj(S) \ar@/^1pc/[l]^{C\otimes_S-} }$$ provides an equivalence of categories, and if $_RC$ is self-co-small then 
$$\xymatrixcolsep{5pc}\xymatrix{ \Prod_R(C^+)\;\ar@/^1pc/[r]^{-\otimes_RC}  &\inj(S) \ar@/^1pc/[l]^{\Hom_S(C,-)} }$$ is an equivalence of categories.
\end{lem}


\begin{cor}\label{special forms} The following assertions hold:
\begin{enumerate}
\item (\cite[Proposition 3.1]{BDGO21}) If $C$ is self-small, then $\Add_R(C)=C\otimes_S\proj(S).$
\item If $_RC$ is self-co-small, then $\Prod_R(C^+)=\Hom_S(C,\inj(S)).$
\end{enumerate}
\end{cor}

\begin{lem}\label{proj are aus-inj are bas} Let $C$ be an $R$-module. The following assertions hold:
\begin{enumerate}
\item $\proj(S)\subseteq \A_C(S)$ if and only if $\text{$_RC$ is $\Sigma$-self-orthogonal and self-small.}$ In this case $\Add_R(C)\subseteq \B_C(R)$.
\item $\inj(S)\subseteq \B_C(S)$ if and only if $\text{$_RC$ is $\prod$-$\Tor$-orthogonal and self-cosmall.}$ In this case $\Prod_R(C^+)\subseteq \A_C(R).$
\item If $_RC$ has a degreewise finite projective resolution, then $\flat (S)\subseteq \A_C(S)$ if and only if $\text{$_RC$ is $\prod$-$\Tor$-orthogonal.}$ In this case $\F(R)\subseteq \B_C(R)$.
\end{enumerate}
\end{lem}
\proof {\it 1.} Follows by \cite[Proposition 5.4(1)]{BGO16} and the fact $\Add_R(C)=C\otimes_S\proj(S)$.

{\it 2.} By the dual argument to that of \cite[Proposition 5.4(1)]{BGO16} and the item {\it 1}.

{\it 3.} Using Proposition  \ref{Gc-proj is Gc-flat}, the "if" part follows by  \cite[Proposition 5.2]{BDGO21} and the "only if" part follows by item {\it 1} since $\proj(S)\subseteq\flat(S)$. \cqfd

\begin{rem}
\begin{enumerate}
\item By Lemma \ref{proj are aus-inj are bas}, the adjoint pair $(C\otimes_S-, \Hom_R(C,-))$ is left semidualizing (in the sense of \cite[Definition 2.1]{GLT15}) if and only if $_RC$ is $\Sigma$-self-orthogonal and self-small. 
\item When $C$ is considered  as a right $S$-module, there is a version for each definition and result presented in this paper.  For example, If $R=\End_S(C)$, then we have the following equalities  $\Add_S(C)=\proj(R)\otimes_RC$ and $\Prod_S(C^+)=\Hom_R(C,\inj(R))$ when $C_S$ is self-small and self-co-small, respectively.
\end{enumerate}
\end{rem}

\begin{lem}\label{Fc-flat cover & Ic-inj-preenvelope} Let $_RC$ be self-co-small and $\prod$-$\Tor$-orthogonal.
\begin{enumerate}
\item Any right $R$-module in $\A_C(R)$ has a monic $\I$-injective preenvelope with cokernel in $\A_C(R)$.
\item Assume that $_RC$ has a degreewise finite projective resolution.  Any $R$-module in $\B_C(R)$ has an epic $\F$-flat cover with kernel in $\B_C(R)$.	
\end{enumerate}
\end{lem}
\proof {\it 1.} By Lemma \ref{proj are aus-inj are bas}(2) we have $\Prod_R(C^+)\subseteq \A_C(R)$, so dual arguments to those of \cite[Proposition 3.8]{BDGO21}, following similar steps to those of \cite[Proposition 3.9]{BDGO21}, proves the statement.

{\it 2.} Let $M\in\B_C(R)$. Then, there exists by \cite[Proposition 3.7(1)]{BEGO22} an $\F$-flat cover $\gamma:L\to M$ which is epic by \cite[Proposition 3.8]{BDGO21}. It remains to show that $\Ker(\gamma)$ is in $\B_C(R)$. But this can be shown exactly as in \cite[Proposition 3.8]{BDGO21}, using the fact that $\F(R)\subseteq \B_C(R)$ by Lemma \ref{proj are aus-inj are bas}(3).\cqfd

\section{Relative cotorsion modules}

In this section we introduce some classes of relative cotorsion modules. Besides their links with other known classes of modules (cotorsion, flat, $\F$-flat, etc.), we are interested in discovering the main homological properties of these new classes. 

Recall that a module $M$ is cotorsion if $\Ext^1(F,M)=0$ for every flat module $F$, equivalently, if $\Ext^i(F,M)=0$ for every flat module $F$ and every $i\geq 1$ (\cite[Definiton 5.3.22]{EJ00}). In the following definition we extend the concept of cotorsion modules to our relative setting. 

\begin{defn} Given an $R$-module $C$ and an integer $n\geq 1$, an $R$-module $M$ is said to be: 
\begin{itemize}
\item $n$-$\C$-cotorsion if $\Ext_R^i(N,M)=0$ for all $\F$-flat modules $N$ and all integer numbers $i$ such that $1\leq i\leq n$. 
\item $\C$-cotorsion if it is $1$-$\C$-cotorsion.
\item Strongly $\C$-cotorsion if it is $n$-$\C$-cotorsion for every $n\geq 1$.
\end{itemize}

The classes of all $\C$-cotorsion modules, of all $n$-$\C$-cotorsion modules and of all strongly $\C$-cotorsion modules will be denoted as $\C(R)$, $\C^n(R)$ and $\SC(R)$, respectively.
\end{defn}

\begin{rems}
\begin{enumerate}
\item When $R$ is a commutative noetherian ring  and $_RC$ is a semidualizing $R$-module, strongly $\C$-cotorsion modules coincide with the $C$-cotorsion modules defined in \cite{SSW11} and the strongly $C$-cotorsion modules defined in \cite{CC16}.
\item Given an integer $n\geq 1$, every $(n+1)$-$\C$-cotorsion $R$-module is $n$-$\C$-cotorsion. Then, we have the following ascending sequence: $$\SC(R) \subseteq \cdots\subseteq \C^{n+1}(R) \subseteq \C^n(R)  \subseteq \cdots \subseteq \C(R),$$ where $\SC(R)$ can be written as $\SC(R)=\cap_{n\geq 1}\C^n(R)$.
\end{enumerate}
\end{rems}

\begin{exmps}
\begin{enumerate}
\item Denote by $\cot(R)$ the class of all cotorsion $R$-modules. If $_RC$ is a flat generator, then  $\cot(R)=\C(R)=\SC(R)$.
\item Every injective module is (strongly) $\C$-cotorsion.
\item Assume that $_RC$ is $\prod$-$\Tor$-orthogonal. Given any $\F$-cover $\varphi:F\to M$ (which exists by \cite[Proposition 3.7]{BEGO22}), $\Ker (\varphi)$ is  $\C$-cotorsion (\cite[Proposition 2.2.7]{Xu96}).
\item Recall that a module $M$ is called Gorenstein cotorsion if $\Ext_R^1(G,M)=0$ for every Gorenstein flat module $G$. If $_RC$ is flat,  then $\F(R)\subseteq \flat(R)\subseteq \Gflat(R)$. Hence, both cotorsion and Gorenstein cotorsion modules are $\C$-cotorsion.
\item  Assume that every $\F$-flat $R$-module has finite injective dimension. Then, $\inj(R)^{\perp_\infty}\subseteq \SC(R)$. In particular, every Gorenstein injective $R$-module is (strongly) $\C$-cotorsion.  
\end{enumerate}
\end{exmps}

It is unknown whether or not the class of $\C$-cotorsion $R$-modules is closed under cokernels of monomorphisms. However, when this happens, all the classes of relative cotorsion $R$-modules introduced above coincide.

\begin{prop} \label{prop of SC-cot}
The following statements hold for any $R$-module $C$.
\begin{enumerate}
\item The class $\SC(R)$ is closed under cokernels of monomorphisms.
\item Let $n\geq 1$ be any integer number. The following assertions are equivalent:
\begin{enumerate}
\item The class $\C^n (R)$ is closed under cokernels of monomorphisms.
\item For every integer number $k\geq n$, every $k$-$\C$-cotorsion module is  $(k+1)$-$\C$-cotorsion.
\item For every integer number $k\geq n$, the class $\C^k (R)$ is closed under cokernels of monomorphisms.
\item Every $n$-$\C$-cotorsion module is strongly $\C$-cotorsion.
\end{enumerate}

In this case, $ \SC(R)= C^{k}(R)$ for every $k\geq n$.
\end{enumerate}
\end{prop}
\proof  (1) Straightforward.

(2) {\it (a)} $\Rightarrow$ {\it (b)} Let $X$ be $n$-$\C$-cotorsion and consider a short exact sequence of $R$-modules $$0\to X\to I\to L\to 0$$ with $I$ injective. By hypothesis, $L$ is $n$-$\C$-cotorsion, so the exact sequence $$0=\Ext_R^n(F,L)\to \Ext_R^{n+1}(F,X)\to \Ext_R^{n+1}(F,I)=0$$ shows that $\Ext_R^{n+1}(F,X)=0$ for every $F\in\F(R)$, that is $X$ i $(n+1)$-$\C$-cotorsion.

But then $L$ is also $(n+1)$-$\C$-cotorsion, and repeating the argument we see that $X$ is indeed $k$-$\C$-cotorsion for every $k\geq n$.

\medskip
{\it (b)} $\Rightarrow$ {\it (c)} Let $0\to K\to M\to N\to 0$ be any exact sequence with $K, M\in \C^k(R)$, choose any $F\in \F (R)$ and consider the induced exact sequence $$\Ext^k(F,M)\to \Ext^k(F,N)\to \Ext^{k+1}(F,K).$$

If $k<n$ then $\Ext^k(F,M)=\Ext^{k+1}(F,K)=0$ and then $\Ext^{k+1}(F,N)=0$.

If $k\geq n$ then by repeatedly applying (2) we get that $M\in \C^k(R)$ and $K\in \C^{k+1} (R)$, so again $\Ext^k(F,M)=\Ext^{k+1}(F,K)=0$ and then $\Ext^{k+1}(F,N)=0$.

Thus, in any case $N\in \C^k(R)$.

\medskip
Now, {\it (c)} $\Rightarrow$ {\it (a)} and {\it (b)} $\Leftrightarrow$ {\it (d)} are clear. \cqfd 

In light of Examples \ref{exmps of Fc-flat}, it is natural to wonder whether there is a relation between $\C$-cotorsion $R$-modules and cotorsion $S$-modules. The following result gives a useful relation.

\begin{prop} \label{Coto and C-cot}
Assume that $C$ is a finitely presented $R$-module and let $n\geq 1$, be an integer. An $R$-module $M$ is $n$-$\C$-cotorsion if and only if $M\in C^{\perp_n}$ and $\Hom_R(C,M)$ is a cotorsion $S$-module. Consequently, 
	$\B_C(R)\cap\C^n(R)=C\otimes_S \left( \A_C(S)\cap\cot(S)\right) $
\end{prop}
\proof $(\Rightarrow)$ Assume that  $M$ is $n$-$\C$-cotorsion. Clearly  $M\in C^{\perp_n}$, since $_RC$ is $\F$-flat. We prove now that $\Hom_R(C,M)$ is a cotorsion $S$-module. 

Let $F$ be any flat $S$-module and consider an exact sequence of $S$-modules $$0\to K\to P\to F\to 0$$ with $P$ projective.  Since the last sequence is pure we have that $K$ is flat and that the sequence $$0\to C\otimes_S K\to C\otimes_S P\to C\otimes_S F\to 0$$ is exact. Then, (Example \ref{exmps of Fc-flat}) $C\otimes_SF$ is $\F$-flat, so $\Ext^1(C\otimes_S F,M)=0$.

Thus, we have the following commutative diagram with exact rows: $$\xymatrixcolsep{1pc}\xymatrix{ \Hom_R(C\otimes_S P, M)\ar[r] \ar[d]^{\cong} &  \Hom_R(C\otimes_S K, M) \ar[r] \ar[d]^{\cong} & 0 \\ \Hom_S(P,\Hom_R(C,M)) \ar[r] & \Hom_S(K,\Hom_R(C,M))\ar[r]   & \Ext^1_S(F,\Hom_R(C,M)) \ar[r] & 0}$$

Therefore, $\Ext^1_S(F,\Hom_R(C,M))=0$ and hence $\Hom_R(C,M)$ is cotorsion.

\medskip
$(\Leftarrow)$ Consider an exact sequence of $R$-modules $0\to M\to I\to L\to 0$, where $I$ is injective. Since $\Ext_R^1(C,M)=0$, the induced sequence $$0\to \Hom_R(C,M)\to \Hom_R(C,I)\to \Hom_R(C,L)\to 0$$ is exact. 

Let $C\otimes_S F$ be any $\F$-flat $R$-module and let's proceed by induction on $n$.

For $n=1$: By the the implication $(\Rightarrow) $, $\Hom_R(C,I)$ is cotorsion since $I$ is $\C$-cotorsion. Then, $\Ext^1_R(F,\Hom_R(C,I))=0$. Also, we have by hypothesis that $\Ext_R^1(F,\Hom_R(C,M))=0$, so the commutative diagram
$$\xymatrixcolsep{1pc}\xymatrix{ \Hom_R(C\otimes_S F, I)\ar[r] \ar[d]^{\cong} &  \Hom_R(C\otimes_S F, L) \ar[r] \ar[d]^{\cong}   & \Ext^1_R(C\otimes_S F, M) \ar[r]  & 0 \\ \Hom_S(F,\Hom_R(C,I)) \ar[r] & \Hom_S(F,\Hom_R(C,L))\ar[r] & 0}$$ has exact rows and then $\Ext^1_R(V, M)=\Ext^1_R(C\otimes_S F, M)=0$, that is, $M$ is $\C$-cotorsion.

Assume now that $n>1$ and that $M\in C^{\perp_{n+1}}$. By induction we know that $M$ is $n$-$\C$-cotorsion, so we only need to prove that $\Ext_R^{n+1}(C\otimes_S F,M)=0$.

Since $_RI$ is injective (and then $(n+1)$-$\C$-cotorsion), $\Hom_R(C,I)$ is a cotorsion $S$-module and then the exact sequences $$0=\Ext^k_R(C,I)\to \Ext^k_R(C,L)\to \Ext^{k+1}_R(C,M)=0$$ and $$0=\Ext_S^k(F,\Hom_R(C,I))\to \Ext_S^k(F,\Hom_R(C,L))\to  \Ext^{k+1}_S(F,\Hom_R(C,M))=0 $$ show that $\Ext_S^k(F,\Hom_R(C,L))=0= \Ext^k_R(C,L)$ for every $k=1,\dots ,n$. Using induction again, we get that $L$ is $n$-$\C$-cotorsion, so the exact sequence $$0=\Ext^n_R(C\otimes_S F,L)\to \Ext^{n+1}_R(C\otimes_S F,M)\to \Ext^{n+1}_R(C\otimes_S F,I)=0,$$ gives that $\Ext_R^{n+1}(C\otimes_S F,M)=0$ as desired.

We now prove the equality $\B_C(R)\cap\C^n(R)=C\otimes_S \left( \A_C(S)\cap\cot(S)\right)$.

Let $M$ be an $R$-module. If $M\in\B_C(R)\cap\C^n(R)$, then $M=C\otimes_S F$ with $_SF\in\A_C(S)$.  Moreover, since $M\in\C^n(R)$, $F\cong \Hom_S(C,C\otimes_SF)=\Hom_S(C,M)$ is a cotorsion $S$-module by the above equivalence. Hence, $M=C\otimes_SF\in C\otimes_S (\A_C(S)\cap\cot(S))$.

For the other inclusion, assume that $M=C\otimes_S F$ with $_SF\in\A_C(S)\cap {\mathcal C}(S)$. Clearly, $M$ is in $\B_C(R)$. On the other hand, we have that $\Hom_R(C,M)\cong F$ is a cotorsion $S$-module and  $\Ext^i(C,M)=0$ for every $i\geq 1$, since $M\in\B_C(R)$. Using again the equivalence proved above, we get that $M\in \C^n(R)$. \cqfd


\begin{cor} \label{Coto and strong C-cot}
Assume that $_RC$ is finitely presented. An $R$-module $M$ is strongly $\C$-cotorsion if and only if $M\in C^{\perp_\infty}$ and $\Hom_R(C,M)$ is a cotorsion $S$-module. Consequently, 
$\B_C(R)\cap\SC(R)=C\otimes_S \left( \A_C(S)\cap\cot(S)\right)=\B_C(R)\cap\C(R)$.
\end{cor}

Given a regular cardinal number $\kappa$, following \cite[Ddefinition 2.1]{EL02} and \cite[Definition 3.6]{LDY11}, a class $\mathcal{A}$ of $R$-modules or complexes of $R$-modules is a $\kappa$-Kaplansky class if for every object $M\in\mathcal{A}$ and for each $x\in M$, there exists a subobject $N$ of $M$ that contains $x$ with the property that $|N|\leq \kappa$ and both $N$ and $M/N$ are in $\mathcal{A}$. We say that $\mathcal{A}$ is a Kaplansky class if it is a $\kappa$-Kaplansky class for some regular cardinal number $\kappa$.

Two proofs were given by  Bican, El Bashir and Enochs in \cite{BEE01}, showing that the class of flat modules forms the left side of a perfect cotorsion pair ($\flat(R)$,$\cot(R)$). In that paper, the authors solved what was known as the flat cover conjecture. In the following result, whose proof is inspired by that given by Enochs, we prove a relative version of that conjecture. This is the first ingredient to get our first Hovey triple in this paper.

\begin{thm}\label{cot pair-relative flat}
Let $C$ be $\prod$-$\Tor$-orthogonal $R$-module. The following statements hold.
\begin{enumerate}
\item  $\left(^\perp\C(R),\C(R)\right)  $ is a complete cotorsion pair cogenerated by a set. Moreover, every $R$-module has a $\C$-cotorsion envelope.
\item $\left(^\perp\SC(R),\SC(R)\right)  $ is a complete  hereditary cotorsion pair cogenerated by a set.
\item  The following assertions are equivalent:
\begin{enumerate}
\item $\left( \F(R),\F(R)^\perp\right)  $ is a perfect cotorsion pair cogenerated by a set.
\item Every $R$-module has a special $\F$-flat precover.
\item Every flat $R$-module is $\F$-flat.
\item $R$ is an $\F$-flat $R$-module.
\end{enumerate}
\end{enumerate}
\end{thm}

\proof  
{\it 1.} Let $F\in\F(R)$. If $\kappa\geq |R|$ then for each $x\in F$  (take $N=Rx$ and $f:N\hookrightarrow M$ in \cite[Lemma 5.3.12]{EJ00}) there is a pure submodule $F_0\leq F$ with $x\in F_0$ such that $|F_0|\leq \kappa$. Since the class $\F(R)$ is closed under pure submodules and pure quotients by \cite[Proposition 3.7(2)]{BEGO22}, $F_0,F/F_0\in \F(R)$.

The latter means that $\F(R)$ is a $\kappa$-Kaplansky class. So, using transfinite induction, we can write $F$ as the direct union of a continuous chain $(F_\alpha)_{\alpha\leq \lambda}$ of pure submodules of $F$ such that $F_0,F_{\alpha+1}/F_\alpha\in\F(R)$ and  $|F_0|,\ |F_{\alpha+1}/F_\alpha|\leq \kappa$ whenever $\alpha+1<\lambda$. 

Let $X$ be the direct sum of all  representatives of $\F(R)$ with cardinality at most $\kappa$. Clearly, $\C(R)\subseteq X^\perp$, and \cite[Theorem 7.3.4]{EJ00} gives that $X^{\perp}\subseteq \F(R)^{\perp}=\C(R)$. Thus, $\C(R)=X^\perp$ and hence  $\left(^\perp\C(R),\C(R)\right)$ is a complete cotorsion pair by \cite[Theorem 7.4.1]{EJ00}. 

The last statement follows from \cite[Theorem 2.8]{EL02}.

{\it 2.} If we prove that $\SC(R)=M^\perp$ for some $R$-module $M$, this assertion will follow from Proposition \ref{prop of SC-cot} and \cite[Theorem 7.4.1]{EJ00}.

By the proof of item 1{\it 1} we have $\C(R)=X^{\perp}$. We claim that $\SC(R)= X^{\perp_\infty}$.

Clearly, $\SC(R)\subseteq X^{\perp_\infty}$. Conversely,  take $N\in X^{\perp_\infty}$ and let $0\to N\to I\to L\to 0$ be a short exact sequence of $R$-modules where $I$ is injective. Note that $L\in X^{\perp_\infty}\subseteq X^\perp=\C(R)$, so by the long exact sequence we get that $$0=\Ext^1_R(F,L)\to \Ext^2_R(F,N)\to \Ext^2_R(F,I)=0$$ for every $F\in\F(R)$. Hence, $\Ext^{2}_R(F,N)=0$. Repeating this process, we get that $\Ext^{i}_R(F,N)=0$ for every $i\geq 1$. Therefore, $\SC(R)= X^{\perp_\infty}$. This means by Lemma \ref{gen-cog sets} that $\SC(R)=X^{\perp_\infty}=M^\perp$  for some $R$-module $M$.


{\it 3.} The implications {\it (a)}$\Rightarrow${\it (b)} and {\it (c)}$\Rightarrow${\it (d)} are obvious.

{\it (b)}$\Rightarrow${\it (c)} Let $F$ be a flat $R$-module and consider a special $\F$-flat precover of $F$: $0\to K\to X\to F\to 0.$ Since $F$ is flat, this sequence is pure. But $\F$ is closed under pure quotients so $F$ is $\F$-flat.

{\it (d)}$\Rightarrow${\it (a)} By (1), \cite[Theorem 5.2.3]{EJ00}, and \cite[Proposition 3.7]{BEGO22}, we only need to show that $\F(R)=\,^\perp \C(R)$.

Clearly, $\F(R)\subseteq \,^\perp \C(R)$. Conversely, take $X \in \;^\perp\C(R)$ and consider an $\F$-flat cover $f:F\to X$, which exists by \cite[Proposition 3.7(1)]{BEGO22}. Since the class $\F(R)$ is closed under direct sums and summands and $_RR$ is $\F$-flat, we get that $\proj(R)\subseteq \F(R)$. Hence, the morphism $f$ is surjective and $K=\Ker (f)$ is $\C$-cotorsion by Lemma \ref{Wakamtsu Lemma} (recall that $\F$ is closed under extensions by \cite[Proposition 3.7]{BEGO22}). But since  $X \in\, ^\perp \C(R)$, $X$ is a direct summand of $F\in\F(R)$. Hence, $X\in\F(R)$ and thus $\F(R)=\,^\perp \C(R)$.\cqfd

\section{Relative Gorenstein flat model structures}

In this section we prove that the class of $\GC$-flat modules is the left hand class  of a perfect hereditary cotorsion pair. Consequently, every module has a $\GCF$-cover. Using Hovey's one-to-one correspondence between abelian model structures and Hovey triples, we obtain a model structure on $R$-Mod whose cofibrant objects are precisely the $\GC$-flat modules. But first we introduce and investigate some properties of Gorenstein cotorsion modules with respect to $C$.

\begin{defn}
Let $C$ be an $R$-module. An $R$-module $M$ is said to be $\GC$-cotorsion if $\Ext_R^1(N,M)=0$ for all $\GC$-flat modules $N$.

We use $\GCC(R)$  to denote the class of all $\GC$-cotorsion $R$-modules.	
\end{defn}

\begin{rems}
When $C=R$, $\GCC(R)$ coincides with the class of all Gorenstein cotorsion modules $\Gcot(R)=\Gflat(R)^\perp$ ($\Gflat(R)$ is the class of all Gorenstein flat modules).
\end{rems}

$\GC$-cotorsion modules are both cotorsion and (strongly) $\C$-cotorsion, as the following result shows.

\begin{prop} \label{cot and Gw-cot rela}
Let $C$ be a $\prod$-$\Tor$-orthogonal $R$-module. Then, every $\GC$-cotorsion module is $\C$-cotorsion. Moreover, if $R$ is $\GCF$-closed and $C$ is $w^+$-tilting, then $$\GCC(R)=\cot(R)\cap\SC(R)\cap \mathcal{H}_C(R),$$ where  $\mathcal{H}_C(R)$ is the class of $R$-modules $M$ such that the complex $\Hom_R(\textbf{X},M)$ is exact for all complete $\F$-flat sequences $\textbf{X}$.
\end{prop}
\proof The first statement follows from the inclusion $\F(R)\subseteq \GCF(R)$ by \cite[Corollary 4.9]{BEGO22}. Now we prove the equality.

$(\subseteq)$ We know that ${\mathcal F}(R)\subseteq \GCF(R)$ by \cite[Proposition 4.17]{BEGO22}, so we have $$\GCC(R)= \GCF(R)^{\perp}\subseteq {\mathcal F}(R)^{\perp}= {\mathcal C}(R).$$

On the other hand, we also know that $\F(R)\subseteq \GCF(R)$, so $$\GCF(R)^{\perp_{\infty}} \subseteq \F(R)^{\perp_{\infty}}=\SC(R).$$ Thus, if $M$ is any $G_C$-cotorsion module, to prove that $M$ is strongly $\C$-cotorsion we only need to prove that $\Ext^{\geq 2}_R (G,M)=0$ for every $G\in \GCF(R)$. Choose then any such $G$ and any projective resolution $$\cdots \to P_1\to P_0\to G\to 0$$ of $G$. If we call $K_i=\Ker(P_{i-1}\to P_{i-2})\ \forall i\geq 1$, we see by \cite[Theorem 4.21]{BEGO22} that $K_i\in\GCF(R)\ \forall i$ ($G$ and all the $P_i$'s are $G_C$-flat). Then, we have, for every $i\geq 2$, $$\Ext_R^i(G,M)\cong \Ext^1_R(K_{i-1},M)=0.$$

Finally, choose any complete $\F$-flat resolution $$\textbf{X}:\cdots \to F_1\to F_0\to C_{-1}\to C_{-2}\to \cdots$$ and let us see that $\Hom_R(\textbf{X},M)$ is exact for any $G_C$-cotorsion module $M$.

By \cite[Corollaries 4.9 and 4.19]{BEGO22} every image $I_i=\Im(F_{i+1}\to F_i)$ and every kernel $K_j=\Ker( C_{j}\to C_{j-1})$ are $\GC$-flat, so $\Ext_R^1(I_i,M)=0=\Ext_R^1(K_j,M)$ for all $i\geq 0$ and all $j\leq -1$, which implies that $\Hom_R(\textbf{X},M)$ is exact.

\medskip
$(\supseteq)$ Let $M\in \cot(R)\cap\SC(R)\cap \mathcal{H}_C(R)$ and $N$ be $\GC$-flat. Then, there exists a complete $\F$-flat resolution $\textbf{X}$ of $N$ as above. Consider the short exact sequence $0\to I_0\to F_0\to N\to 0.$ Since this sequence is $\Hom_R(-,M)$-exact by the hypotheses and $M$ is cotorsion, the exactness of the sequence $$0\to \Hom_R(N,M)\to \Hom_R(F_0,M)\to \Hom_R(I_0,M)\to \Ext_R^1(N,M)\to 0$$ shows that $\Ext_R^1(N,M)=0$. Thus, $M$ is $\GC$-cotorsion. \cqfd

Now we are in position to investigate when the classes $\GCF(R)$ and $\GCC(R)$ are  covering and enveloping, respectively. In fact, we will prove more than that: we will characterize exactly when $(\GCF(R),\GCC(R))$ is a perfect and hereditary cotorsion pair. We will use the following lemma. 

\begin{lem} \label{Kaplansky}
Assume that $R$ is $\GCF$-closed and  $C$ is $w^+$-tilting. Then, the class of $\GC$-flat $R$-modules is a Kaplansky class.
\end{lem}
\proof Let $\tilde{\mathcal{A}}$ be the class of all $\left(\Prod_R(C^+)\otimes_R-\right)$-exact exact complexes of $R$-modules with components in $\mathcal{A}:=\flat(R)\cup \F(R)$. Since  $\F(R)$ and $\flat(R)$ are closed under direct summands, pure submodules and pure quotients (see \cite[Proposition 3.7 (2)]{BEGO22}), so is the class $\mathcal{A}$. On the other hand, as in the proof of \cite[Theorem 3.7]{EG19}, we get that the class $\tilde{\mathcal{A}}$ is closed under pure subcomplexes and pure quotients (in the sense of \cite[Definition 3.1]{EG19}). Using now  \cite[Proposition 3.4]{EG19}, we get that the class $\tilde{\mathcal{A}}$ is a Kaplansky class. But since $\GCF(R)$ is the class of $0$-syzygies of exact sequences in $\widetilde{\mathcal{A}}$ by \cite[Corollary 5.2]{BEGO22}, we get that it is also a Kaplansky class, as desired. \cqfd 

\begin{cor}
Assume that $R$ is $\GCF$-closed and  $C$ is $w^+$-tilting. Then,  $\GCF(R)$ is preenveloping if and only if it is closed under direct products. 
\end{cor}
\proof Follows from \cite[Theorem 2.5]{EL02} since the class $\GCF$ is closed under direct limits by \cite[Proposition 4.15]{BEGO22}. \cqfd

\begin{thm}\label{(GwF,GwC) cot pair}
Let $C$ be a $\prod$-$\Tor$-orthogonal $R$-module. The following assertions are equivalent:
\begin{enumerate}
\item $(\GCF(R),\GCC(R))$ is a perfect hereditary cotorsion pair cogenerated by a set.
\item $R$ is $\GCF$-closed and  $C$ is $w^+$-tilting.
\end{enumerate}

In this case, $\GCF(R)$ is covering and  $\GCC(R)$ is enveloping.
\end{thm}
\proof $2.\Rightarrow 1.$ Using Lemma \ref{Kaplansky}, $\GCF(R)$ is a Kaplansky class. By \cite[Propositions 4.15 and 4.17]{BEGO22}, it is closed under direct limits and contains all projective $R$-modules. Therefore, since $\GCF(R)$ is closed under extensions, we get that $(\GCF(R),\GCC(R))$ is a perfect cotorsion pair by \cite[Theorem 2.9]{EL02}. Moreover, since the class $\GCF(R)$ is projectively resolving by \cite[Theorem 4.21]{BEGO22}, our pair is hereditary. Finally, as in the proof of Theorem \ref{cot pair-relative flat}(1), our pair is cogenerated by a set.

$1.\Rightarrow 2.$ By hypothesis, $^\perp \GCC(R)= \GCF(R)$, so the class $\GCF(R)$ is closed under extensions and since $R\in \GCF(R)$,  we get that $C$ is $w^+$-tilting by \cite[Proposition 4.17]{BEGO22}. \cqfd

\begin{lem}\label{exmp of GwF-closed}
Assume that $_RC$ is finitely presented. Then, $\F(R)$ is closed under direct products if and only if  $S$ is right coherent and $C_S$ is finitely presented.
	
In this case, if $C$ is $\prod$-$\Tor$-orthogonal, then $R$  is $\GCF$-closed.
\end{lem}
\proof $(\Rightarrow)$ Since $_RC$ is finitely presented, Example \ref{exmps of Fc-flat} says that $\F(R)=C\otimes_S \flat (S)$. Then, for any $\F$-flat $R$-module $C\otimes_S F$ we have that the canonical homomorphism $\tau_{C\otimes_S F}:C\otimes_S\Hom_R(C,C\otimes_S F)\to C\otimes_S F$ is an isomorphism by \cite[Theorem 3.2.14]{EJ00}. In particular, since $C^I\in \F(R)$ for any set $I$ by the hypotheses, $\tau_{C^I}$ is always an isomorphism. But $\tau_{C^I}$ is nothing but the composition of the canonical morphisms $$C\otimes_S \Hom_R(C,C^I)\stackrel{\beta}{\to} C\otimes_S\Hom_R(C,C)^I=C\otimes_SS^I\stackrel{\alpha}{\to} (C\otimes_S S)^I=C^I$$ and $\beta$ is an isomorphism. Thus, $\alpha$ is also an isomorphism for any set $I$ and then $C_S$ is finitely presented.

Now, since $C^I\in \F(R)$, there must exist some $F\in\flat (S)$ such that $C^I=C\otimes_S F$. Then,we can use \cite[Theorem 3.2.14]{EJ00} again to get the following chain of isomorphisms: $$_SF\cong S\otimes_S F=\Hom_R(C,C)\otimes_S F\cong \Hom_R(C,C\otimes_S F)=$$ $$=\Hom_R(C,C^I)\cong \Hom_R(C,C)^I=S^I.$$

Therefore, $_SS^I$ is flat for any set $I$ and then $S$ is right coherent.

\medskip
$(\Leftarrow)$ By Example \ref{exmps of Fc-flat} we know that $\F(R)=C\otimes_S \flat (S)$. Thus, choose any family $\{ F_i;\ i\in I\}$ of flat $S$-modules and let us prove that $\prod_{i\in I} C\otimes_S F_i\in \F(R)$. But $C_S$ is finitely presented so $\prod_{i\in I} C\otimes_S F_i\cong C\otimes_S(\prod_i F_i)$, and $\prod_i F_i$ is flat because $S$ is right coherent.

Finally, the fact that $R$ is $\GCF$-closed follows from \cite[Corollary 4.13]{BEGO22}.\cqfd

\begin{cor}
Let $C$ be a semidualizing  $(R,S)$-bimodule such that $S$ is right coherent. The following assertions hold: 
\begin{enumerate}
\item $(\GCF(R),\GCC(R))$ is a perfect and hereditary cotorsion pair.
\item $\GCF(R)$ is special precovering.
\item $\GCC(R)$ is special preenveloping.
\end{enumerate}
\end{cor} 
\proof Follows from Theorem \ref{(GwF,GwC) cot pair}, Lemma \ref{Wakamtsu Lemma} and Lemma \ref{exmp of GwF-closed}. \cqfd

\begin{cor}
Assume that $R$ is $\GCF$-closed and that $C$ is a $w^+$-tilting $R$-module. Then, the following assertions are equivalent:
\begin{enumerate}
\item Every $R$-module is $\GC$-cotorsion.
\item Every $\GC$-flat $R$-module is $\GC$-cotorsion.
\item Every $\GC$-flat $R$-module is projective.
\item $R$ is left perfect and every $\GC$-flat $R$-module is flat.
\item $R$ is left perfect and every cotorsion $R$-module is $\GC$-cotorsion.
\end{enumerate}
\end{cor}
\proof  $1. \Rightarrow 2.$ Clear.

$2. \Rightarrow 3.$ Let $M$ be a $\GC$-flat $R$-module and consider an exact sequence $0\to K\to P\to M\to 0$ where $P$ is projective. Since $P$ is $\GC$-flat by \cite[Proposition 4.17]{BEGO22} and the class $\GCF(R)$ is closed under kernels of epimorphisms by \cite[Proposition 4.15]{BEGO22}, $K$ is $\GC$-flat and then $\GC$-cotorsion by hypothesis. Therefore, the above sequence is split and $M$ is projective.


$3. \Rightarrow  4.$ Every flat module is $\GC$-flat by \cite[Proposition 4.17(4)]{BEGO22}, so it is projective by hypothesis. Then, $R$ is left perfect and then every $\GC$-flat module is flat 

$4.\Rightarrow 5.$ Since both $(\F(R),{\mathcal C}(R))$ and $\GCF(R),\GCC(R))$ are cotorsion pairs, we have ${\mathcal C}=\flat (R)^{\perp}\subseteq \GCF(R)^{\perp}=\GCC(R)$, where the middle inclusion holds because $\GCF(R)\subseteq \flat (R)$.

$5. \Rightarrow  1.$ Since $R$ is perfect, every $R$-module is cotorsion by \cite[Proposition 3.3.1]{Xu96}. Hence, every $R$-module is $\GC$-cotorsion. \cqfd
 
The following result is inspired by an argument due to Estrada, Iacob and P\'{e}rez in \cite[Proposition 4.1]{EIP20}.

\begin{prop}\label{(Fc,Cc) core}
Assume that $R$ is $\GCF$-closed and that $C$ is a $w^+$-tilting $R$-module. Then, $$\GCF(R)\cap\GCC(R)=\F(R)\cap \C(R).$$
\end{prop} 
\proof  $(\subseteq)$  Assume that $M\in \GCF(R)\cap \GCC(R)$. Then,  $M\in \C(R)$ by Proposition \ref{cot and Gw-cot rela}. Moreover, since $M$ is $\GC$-flat, there exists by \cite[Proposition 4.10]{BEGO22} an exact sequence of $R$-modules $0\to M\to V\to G\to 0$ where $V$ is $\F$-flat and $G$ is $\GC$-flat. This short exact sequence splits since $G$ is $\GC$-flat and $M$ is $\GC$-cotorsion. Hence, $M\in \F(R)$.

$(\supseteq)$ Assume now that $M\in \F(R)\cap \C(R)$. Clearly $M\in \GCF(R)$ (recall that $\F (R)\subseteq \GCF(R)$ by \cite[Corollary 4.9]{BEGO22}). let us prove that $M\in \GCC(R)$.

By Theorem \ref{(GwF,GwC) cot pair}, $M$ has a special $\GCC(R)$-preenvelope $$0\to M\to X\to G\to 0.$$ Since $M$ and $G$ are $\GC$-flat modules and $R$ is $\GCF$-closed, $X$ is $\GC$-flat as well and then $X\in \GCF(R)\cap \GCC(R)\subseteq \F(R)\cap \C(R)$ by the first inclusion. Now, since $M$ and $X$ are $\F$-flat and $G$ is $\GC$-flat, $G\in \F(R)$ by \cite[Theorem 7.12]{BEGO22}. Then, the short exact sequence splits. Hence, $M$ is $\GC$-cotorsion.  \cqfd

We are now ready to show that $\GCF(R)$ (resp., $\GCC(R)$) is part of a left (resp., right) weak AB-context. Note  that $\F(R)\cap \C(R)=\F(R)\cap \SC(R)$. The claim concerning the class $\GCF(R)$ in the following result was proved by Sather-Wagstaff, Sharif and White in \cite{SSW11}, when the  ring $R$ is commutative and noetherian and $C$ is a semidualizing $R$-module. Here we use a different and recent approach and also we show this result with significantly lower restrictions on $R$ and $C$.

\begin{thm} \label{Gc-flat weak AB-cont}
If $R$ is $\GCF$-closed and $C$ is a $w^+$-tilting $R$-module, then $$\left(\GCF(R),\res (\widehat{\F(R)\cap \C(R)}) ,\F(R)\cap\C(R)\right)$$ is a  left weak AB-context, and $$\left(\F(R)\cap\C(R),\cores (\widehat{\F(R)\cap \C(R)}) ,\GCC(R)\right)$$ is a right weak AB-context.
\end{thm}
\proof By Lemma \ref{left Frobnius}(1), it suffices to show that $\left(\GCF(R),\F(R)\cap\C(R)\right)$ is a left Frobenius pair. But this follows by Theorem \ref{(GwF,GwC) cot pair}, Proposition \ref{(Fc,Cc) core} and  Lemma \ref{left Frobnius}(2).  The claim concerning the class $\GCC(R)$ follows similarly using Lemma \ref{right Frobnius}.\cqfd

An immediate consequence of Theorem \ref{Gc-flat weak AB-cont} is that the modules in $\res(\widehat{\GCF(R)})$ (resp., $\cores(\widehat{\GCC(R)})$) have some special approximations in the sense of Auslander and Buchweitz \cite{AB89}. The following result follows from \cite[Theorem 1.1]{AB89} and its dual. For more interesting consequences of Theorem \ref{Gc-flat weak AB-cont} we refer to \cite[Theorem 1.12.10]{Has20}.

\begin{cor} \label{Gc-flat approx}
Assume that $R$ is $\GCF$-closed and $_RC$ is $w^+$-tilting.
\begin{enumerate}
\item For every $M\in \res(\widehat{\GCF(R)})$ there exist exact sequences $$0\to Y_M\to X_M\to M\to 0 \text{ and } 0\to M\to Y^M\to X^M\to 0$$ with  $X_M,X^M\in \GCF(R)$ and $Y_M,Y^M\in \res (\widehat{\F(R)\cap \C(R)})$.
\item For every very $N\in \cores(\widehat{\GCC(R)})$ there exist exact sequences $$0\to Y^N\to X^N\to N\to 0 \text{ and } 0\to N\to Y_N\to X_N\to 0$$ with $Y_N,Y^N\in \GCC(R)$ and $X_N,X^N\in \cores (\widehat{\F(R)\cap \C(R)})$.
\end{enumerate}
\end{cor}

\begin{prop}\label{2nd descr of the core of (GcF,GcC)}
Assume that $C$ is $w$-tilting $R$-module with a degreewise finite projective resolution. Then, $\left(^\perp(\B_C(R)\cap\C(R)) ,\B_C(R)\cap\C(R)\right) $ is a hereditary complete cotorsion pair cogenerated by a set.
\end{prop}
\proof First of all, note that $C$ is $\prod$-$\Tor$-orthogonal by Proposition \ref{Gc-proj is Gc-flat}(1) and $\B_C(R)\cap\C(R)=\B_C(R)\cap\SC(R)$ by Corollary \ref{Coto and strong C-cot}. Now using  \cite[Theorem 3.10.]{BDGO21}, we see that  $\B_C(R)=\X_1^{\perp_\infty}$ for some set $\X_1$. Similarly,  $\SC(R)= \X_2^{\perp_\infty}$  for some set $\X_2$ by Theorem \ref{cot pair-relative flat}(2). Then, $\B(R)\cap\C(R)=\X_1^{\perp_\infty}\cap \X_2^{\perp_\infty}=(\X_1\cup \X_2)^{\perp_\infty}$. Thus, $\B(R)\cap\C(R)=M^\perp$ for some $R$-module $M$ by Lemma \ref{gen-cog sets}. Hence,  our pair is hereditary and cogenerated by a set and then complete by \cite[Theorem 7.4.1]{EJ00}. \cqfd

Under strong conditions we get a different description of the core of the cotorsion pair $(\GCF(R),\GCC(R))$ which is  the last ingredient to get our first model structure.

\begin{prop}
If $R$ is $\GCF$-closed and $C$ is a $w^+$-tilting $R$-module admitting a degreewise finite projective resolution, then $$\GCF(R)\cap\GCC(R)=\;^\perp (\B_C(R)\cap\C(R))\cap (\B_C(R)\cap\C(R))$$	
\end{prop} 
\proof $(\subseteq)$ By Proposition \ref{(Fc,Cc) core},  $ \GCF(R)\cap\GCC(R)\subseteq \GCC(R)\subseteq \F(R)\cap\C(R)$.  But $\F(R)\subseteq \B_C(R)$ by Lemma \ref{proj are aus-inj are bas}(3). Then, $M\in \F(R)\subseteq\B_C(R)$ and hence $ \GCF(R)\cap\GCC(R)\subseteq \B(R)\cap\C(R)$.

On the other hand, we have $\B(R)\cap\C(R)\subseteq \C(R)$, which implies that $^\perp\C(R)\subseteq \,^\perp(\B(R)\cap\C(R))$. But $\F(R)\subseteq \,^\perp\C(R)$. Hence, $\GCF(R)\cap\GCC(R)\subseteq \,^\perp(\B(R)\cap\C(R))$ by Proposition \ref{(Fc,Cc) core}. 

$(\supseteq)$ Conversely, let $M\in\,^\perp(\B(R)\cap\C(R))\cap (\B(R)\cap\C(R))$. By Proposition \ref{(Fc,Cc) core} we only need to show that $M\in\F(R)\cap\C(R)$.

Clearly, $M\in\C(R)$ and since $M\in\B_C(R)$, Lemma \ref{Fc-flat cover & Ic-inj-preenvelope}(2) says that there exists  an epic $\F$-flat cover $\gamma:\F \twoheadrightarrow M$ with  $K:=\Ker(\gamma)\in\B_C(R)$. Wakamatsu Lemma says that $K\in\C(R)$, that is $K\in \B(R)\cap\C(R)$. Since $M\in\,^\perp (\B(R)\cap\C(R))$, the short exact sequence $0\to K\to F\to M\to 0 $ splits. Thus, $M\in\F(R)$, as desired. \cqfd

Modules in $\mathcal{H}_C:=\B_C(R)\cap\C(R)$ will be called $\mathcal{H}_C$-cotorsion modules. They will play the role of fibrant objects in the following model structure. 

\begin{thm}\label{GwF model structure}
Let $R$ be $\GCF$-closed and $C$ be a $w$-tilting $R$-module admitting a degreewise finite projective resolution. Then, there exists a unique hereditary abelian model structure on $R$-{\rm Mod}, called the \textbf{$\GC$-flat model structure}, as follows:

$\bullet$ The cofibrant objects coincide with the $\GC$-flat $R$-modules.

$\bullet$ The fibrant objects coincide with the $\mathcal{H}_C$-cotorsion $R$-modules.

$\bullet$ The class of trivially cofibrant objects coincide with  $^\perp (\B_C(R)\cap\C(R))$.

$\bullet$ The trivially fibrant objects coincide with the $\GC$-cotorsion $R$-modules.
\end{thm}
\proof  It follows from Theorem \ref{cot pair-relative flat} and Theorem \ref{(GwF,GwC) cot pair} that the pairs $$ \left(\GCF(R),\GCC(R) \right) \text{ and }\left(^\perp(\B_C(R)\cap\C(R)),\B_C(R)\cap\C(R) \right)$$ are complete and hereditary cotorsion pairs. By Proposition \ref{2nd descr of the core of (GcF,GcC)}, these cotorsion pairs have the same core. Let us now show that $\GCC(R)\subseteq \B_C(R)\cap\C(R)$.

The inclusion $\GCC(R)\subseteq \C(R)$ holds by Proposition \ref{cot and Gw-cot rela}. On the other hand, as we did in the proof Proposition\ref{Gc-proj is Gc-flat}(2), we can  get an $\add_R(C)$-coresolution of $R$: $$0\to R \stackrel{t_0}{\to} C_0\stackrel{t_1}{\to} C_1\stackrel{t_2}{\to}\cdots$$ which is $\Hom_R(-,C)$-exact (or equivalent, $\Hom_R(-,\Add_R(C))$-exact since each $C_i$ is finitely generated). Then, $\B_C(R)=(C\oplus (\oplus_{i\geq0} \Coker t_i))^\perp$ by \cite[Theorem 3.10]{BDGO21}. 

Now by \cite[Corollary 2.13]{BGO16}, each $\Coker (t_i)$ is $\GC$-projective. But since each $\Coker (t_i)$ has a degreewise finite projective resolution, it is also  a $\GC$-flat module by Proposition \ref{Gc-proj is Gc-flat}. It follows that $C\oplus (\oplus_{i\geq0} \Coker t_i)\in\GCF(R)$ which implies that $\GCC(R)\subseteq (C\oplus (\oplus_{i\geq0} \Coker t_i))^\perp=\B_C(R)$.  Thus, Theorem \ref{Hovey's correpondence} gives the desired Hovey triple.  \cqfd

\begin{rem} 
\begin{enumerate}
\item  Under Hovey's one-to-one correspondence between abelian model structures and Hovey triples \cite[Theorem 2.2]{Hov02}, the $\GC$-flat model structure is described as follows:
		
$\bullet$ A morphism $f$ is a cofibration (trivial cofibration) if and only if it is a monomorphism with $\GC$-flat cokernel ($\Coker (f)\in \;^\perp (\B_C(R)\cap\C(R))$).

$\bullet$ A morphism $g$ is a fibration (trivial fibration) if and only if it is an epimorphism with $\mathcal{H}_C$-cotorsion kernel ($\Ker (g)\in\GCC(R)$).
		
Moreover, this model structure is cofibrantly generated in the sense of \cite[Section 2.1]{Hov99}. This can be seen by \cite[Lemma 6.7 and Corollary 6.8]{Hov02} since the pairs $\left(^\perp(\B_C(R)\cap\C(R)),\B_C(R)\cap\C(R)\right) $ and $\left(\GCF(R),\GCC(R) \right)$  are cogenerated by a set.
\item When $C$ is $\Sigma$-self-orthogonal admitting a degreewise finite projective resolution, $C$ is $w$-tilting and $R$ is $\GCF$-closed whenever the $\GC$-flat model structure exists. 
This can be easily seen by Theorem \ref{(GwF,GwC) cot pair} and Proposition \ref{Gc-proj is Gc-flat} since the $\GC$-flat model structure implies that $\left(\GCF(R),\GCC(R) \right) $ is a complete and hereditary cotorsion pair.
\end{enumerate} 
\end{rem}

Recall that an exact category (in the sense of Quillen) is an additive category with an exact structure, that is, a distinguished class of ker-coker sequences which are called conflations, subject to certain axioms (\cite[Appendix A]{Kel90}). For example, a full additive subcategory of an abelian category that is closed under extensions is an exact category such that conflations are short exact sequences with terms in the subcategory. All the full additive subcategories in this paper will be considered exact with respect to this canonical exact structure. 

Recall that a Frobenius category is an exact category with enough injectives and projectives and such that the projective objects coincide with
the injective objects. Given a Frobenius category $\mathcal{A}$, we can form the stable category $\underline{\A}:=\A/\sim$, which has the same objects as $\mathcal{A}$ and $\Hom_{\A/\sim}(X,Y)=\Hom_{\A}(X,Y)/\sim$, where $f\sim g$ if and only if $f-g$ factors through a projective-injective object.

The main fact about a Frobenius category $\A$ is that the stable category is canonically triangulated and it encodes the corresponding relative homological algebra on $\A$.

The first part of the following result was recently proved by Hu, Geng, Wu and Li in \cite[Theorem 4.3]{HGWL21} when $R$ is a commutative noetherian ring and $C$ is a semidualizing $R$-module. We obtain it here with a different approach and fewer assumptions.

\begin{cor}\label{GcF-homotopy}
Assume that $R$ is $\GCF$-closed and that $C$ is $w$-tilting admitting a degreewise finite projective resolution. Then,  the category $\A_{c,f}:=\GCF(R)\cap\B_C(R)\cap\C(R)$ is a Frobenius category. The projective-injective objects are exactly the objects in  $\W_{c,f}=\F(R)\cap \C(R)$. Moreover, the homotopy category of the $\GC$-flat model structure is triangle equivalent to the stable category $\underline{\GCF(R)\cap\B_C(R)\cap\C(R)}$. 
\end{cor}
\proof Use Theorem  \ref{GwF model structure} and Proposition \ref{(Fc,Cc) core} together with \cite[Proposition 5.2(4) and Lemma 4.7]{Gil11}.  The last assertion follows from \cite[Corollary 5.4]{Gil11}. See also \cite[Proposition 4.2 and Theorem 4.3]{Gil16}. \cqfd

The Gorenstein flat model structure goes back to Gillespie and Hovey \cite[Theorem 3.12]{GH10} when the ring is Iwanaga-Gorenstein. Recently, \v{S}aroch and \v{S}\v{t}ov\'{\i}\u{c}ek proved in \cite{SS20} the existence of this model structure over any arbitrary ring.

\begin{cor}(\textbf{The Gorenstein flat model structure})
For any finitely generated projective generator $C$ of $R$-{\rm Mod}, there exists a unique hereditary abelian model structure on $R$-\Mod where $\GCF(R)=\Gflat(R)$  is the class of cofibrant objects and $\B_C(R)\cap\C(R)=\cot(R)$ is the class of fibrant objects.
	
In this case, the category	$\Gflat(R)\cap \cot(R)$ is a Frobenius category where the  projective-injective objects are exactly the flat-cotorsion $R$-modules. Moreover, the homotopy category of the Gorenstein flat model structure is triangle equivalent to the stable category $\underline{\Gflat(R)\cap\cot(R)}$. 	
\end{cor}
\proof Clearly $\F(R)=\flat(R)$ since $_RC$ is a projective generator. Hence, $\C(R)=\cot(R)$ and $\GCF(R)=\Gflat(R)=\,^\perp\Gcot(R)$ by \cite[Proposition 4.23]{BEGO22} and \cite[Corollary 4.12]{SS20}. This implies that $R$ is $\GCF$-closed. So, by Theorem \ref{GwF model structure} we have the desired model structure.\cqfd

By the above corollary, a projective generator $_RC$ is an example of a case where the $\GC$-flat model structure exists. At the end of this section we will give another example of $C$ being neither projective nor a generator. 

\medskip
Recall from \cite{Ste70} that the FP-injective dimension of an $R$-module $M$, denoted by FP-$\id_R(M)$, is defined to be the smallest non-negative integer $n$ such that, for every finitely presented $R$-module $F$, the equality $\Ext_R^{n+1}(F,M)=0$ holds. In particular,  $M$ is called FP-injective if FP-$\id_R(M)=0$.  We write $\mathcal{FI}(R)$ to denote the class of all FP-injective $R$-modules. Note that $\mathcal{FI}(R)=\inj(R)$ when $R$ is a (left) noetherian ring.

\begin{lem}\label{Fc=FI}Let $C$ be a semidualizing $(R,R)$-bimodule. Then,  $\F(R)=\mathcal{FI}(R)$ if and only if $R$ is left and right coherent, and both $_RC$ and $C_R$ are FP-injective.
\end{lem}
\proof 
$(\Rightarrow)$ Clearly $_RC$ is $FP$-injective. Moreover, $R$ is left coherent by \cite[Theorem 3.2]{Ste70} since $\F(R)$ is closed under direct limits, and $R$ is right coherent by Lemma \ref{exmp of GwF-closed} since $\mathcal{FI}(R)$ is closed under direct products.

Now we prove that $C_R$ is FP-injective. By \cite[Theorem 1]{CS81}, it suffices to show that $C^+$ is a flat $R$-module. We know that $R^+$ is an injective $R$-module, so $R^+\in \mathcal{FI}(R)=\F(R)$ and hence $R^+=C\otimes_R F$ for some flat module $_RF$. Thus, $C^+\cong (R\otimes_R C)^+\cong \Hom_R(_RC,R^+)\cong \Hom_R(_RC,C\otimes_R F)\cong F\in\flat(R)$.

\medskip
$(\Leftarrow) $ Let $X$ be a left $R$-module. If $X\in\F(R)$, then $X^+\oplus Y=(C^+)^I$ for some module $Y$ and some set $I$. But since $_RC$ is FP-injective and $R$ is left coherent, $(C^+)^I$ is flat by \cite[Theorem 1]{CS81} and \cite[Theorem 3.2.24]{EJ00}. Hence, $X^+$ is flat and then $X\in\mathcal{FI}(R)$.

Conversely, assume that $X\in\mathcal{FI}(R)$. Then, $X^{++}\in\inj(R)$ by \cite[Theorem 1]{CS81} which implies that $X^{++}$ is a direct summand of $(R^+)^J$ for some set $J$.

Note also that $(C^+)^J$ is a flat $R$-module since $C_R$ is FP-injective and $R$ is right coherent.  Hence  $X^{++}\cong C\otimes_R \Hom_R(C,X^{++})$ is a direct summand of $$ C\otimes_R \Hom_R(_RC,(R^+)^J)\cong C\otimes_R(C^+)^J\in C\otimes_R\flat(R)=\F(R).$$

Hence $X^{++}\in\F(R)$. But $X$ is a pure  submodule of $X^{++}$  and $\F(R)$ is closed under pure submodules by \cite[Theorem 3.8]{BEGO22}. Therefore, $X\in\F(R)$ and thus $\F(R)=\mathcal{FI}(R)$.\cqfd



Recall (\cite[Definition 2.1]{MD08}) that an $R$-module $M$ is Gorenstein FP-injective if there exists a $\Hom_R(\mathcal{FI}(R),-)$-exact exact sequence of injective $R$-modules $\cdots \to E_1\to E_0\to E_{-1}\to\cdots,$ with $M=\Im(E_0\to E_{-1})$.

\begin{cor}
Let $R$ be left and right coherent with $FP-\id(_RR)<\infty $.  If $C$ is a semidualizing $(R,R)$-bimodule such that $_RC$ and $C_R$ are FP-injective, then there exists a hereditary abelian model structure on $R$-\Mod such that each $R$-module is cofibrant and the fibrant objects are the Gorenstein FP-injective $R$-modules in $\B_C(R)$.
	
	

\end{cor}
\proof Clearly $C$ is $w$-tilting and using Lemma \ref{exmp of GwF-closed} we see that $R$ is $\GCF$-closed. By Lemma \ref{Fc=FI}, $\F(R)=\mathcal{FI}(R)$ and then  $\SC(R)=\mathcal{FI}^{\perp_\infty}$ is the class of Gorenstein FP-injective modules (\cite[Theorem 2.4]{MD08}). 
Consequently, using Corollary \ref{Coto and strong C-cot} we get that $\B_C(R)\cap\C(R)=\B_C(R)\cap \SC(R)$ is the class of Gorenstein FP-injective that are in $\B_C(R)$. Finally, since $\inj(R)\subseteq \F(R)$ and $\Prod_R(C^+)\subseteq \flat(R_R)$, every $R$-module is $\GC$-flat. Therefore, Theorem \ref{GwF model structure} gives the desired model structure. \cqfd

\begin{exmp}[\cite{BEGO21a}, Example 3.5]
Take the quiver $Q: \bullet\to\bullet\to \cdots \to \bullet$ with $n\geq 1$ vertices and $R=kQ$ the path algebra over an algebraic field $k$. By \cite[example 3.5]{BEGO21a} there are two semidualizing $(R,R)$-bimodules $$C_1=R   \text{ and }  C_2=R^*=\Hom_R(R,k).$$ Note that $C_2$ is neither a flat nor a generator module. There are two hereditary abelian model structures on $R$-\Mod which are as follows: $$\mathcal{M}_1=({\rm G}_{\rm C_1}\mathcal{\rm F}(R),\W_1,\B_{C_1}(R)\cap\cot_{C_1}(R))=(\Gproj(R),\Gproj(R)^\perp,R{\rm -Mod}),$$ $$\mathcal{M}_2=({\rm G}_{\rm C_2}\mathcal{\rm F}(R),\W_2,\B_{C_2}(R)\cap\cot_{C_1}(R))=(R{\rm -Mod},\,^\perp\Ginj(R),\Ginj(R))$$ where $\Gproj(R)$ and $\Ginj(R)$ denote the class of Gorenstein projective and injective $R$-modules, respectively.
\end{exmp}

\section{Auslander and Bass model structures}

In this section, we study the existence of  an abelian model structure on $R$-Mod (resp., Mod-$R$) such that $\B_C(R)$ (resp., $\A_C(R)$) is the class of fibrant (resp., cofibrant) objects.

In order to apply Hovey's one-to-one correspondence between abelian model structures and Hovey triples, we first need to know when the Bass (resp., Auslander) class $\B_C(R)$ (resp., $\A_C(R)$) forms the right (resp., left) hand of a complete and hereditary cotorsion pair. This has been recently proved in \cite[Corollary 3.12 and Theorem 3.16(ii)]{BDGO21}.

\begin{prop}\label{Ac and Bc cortorsion pairs}
Let $C$ be an $R$-module. The following assertions hold:
\item \begin{enumerate}
\item If $_RC$ is self-small and w-tilting, then  $(^\perp\B_C(R),\B_C(R))$ is a complete and hereditary cotorsion pair in $R$-{\rm Mod}.
\item  If $_RC$ is self-co-small and  $w^+$-tilting, then $(\A_C(R),\A_C(R)^\perp)$ is a complete and hereditary cotorsion pair in {\rm Mod}-$R$.
\end{enumerate}
\end{prop}
\proof 1. See \cite[Corollary 3.12]{BDGO21}.

2. The proof is almost dual to that  of \cite[Corollary 3.12]{BDGO21}, we only mention a few specific points concerning this dualization. We split the proof into three parts:

$\bullet$ Any $M\in \A_C(R)$ has a monic $\I$-injective preenvelope with cokernel in $\A_C(R)$ by Lemma \ref{Fc-flat cover & Ic-inj-preenvelope}. 

$\bullet$ $(\A_C(R),\A_C(R)^\perp)$ is a complete cotorsion pair: given any $\left( \Prod_R(C^+)\otimes_R-\right)$-exact exact sequence $0\to R\stackrel{t_0}{\to} C_0\stackrel{t_1}{\to} C_1\stackrel{t_1}{\to} \cdots$ with each $C_i$ $\F$-flat and applying the functor $\Hom(-,\mathbb{Q}/\Z)$, we get a $\Hom_R(\Prod_R(C^+),-)$-exact exact sequence $\cdots\to C_1^+\to C_0^+\to R^+\to 0$ where each $C_i^+$ is $\I$-injective. If we let  $K_i=\Coker (t_i)$ with $i\geq 0$ and $K=C\oplus (\oplus_{i\geq 0}K_i)$, then $K^+\cong C^+\oplus (\prod_{i\geq 0}\Ker (t_i)^+)$. But $C^+$ is a w-cotilting right $R$-module so by dualizing the proof of \cite[Theorem 3.10]{BDGO21} we get $\A_C(R)=\;^\perp(K^+)$. Thus, this part follows by \cite[Theroem 7.4.1]{EJ00}.

$\bullet$  The hereditary property  follows  by \cite[Proposition 5.4.(4)]{BGO16}.\cqfd

\begin{prop} \label{Bass-core-cot}
Let $C$ be an $R$-module. The following assertions hold.
\begin{enumerate}
\item If $_RC$ is self-small and $\Sigma$-self-orthogonal, then $^\perp\B_C(R)\cap\B_C(R)=\Add_R(C).$
\item  If $_RC$ is self-co-small and $\prod$-$\Tor$-orthogonal then we have $\A_C(R)\cap\A_C(R)^\perp=\Prod_R(C^+)$.
\end{enumerate}
\end{prop}
\proof 1. Follows by \cite[Proposition 4.9(1)]{BDGO21}  and Lemma \ref{special forms}(2). 

2. Note that $C^+$ is $\prod$-self-orthogonal (\cite[Proposition 3.5]{BDGO21}). So, this equality follows by Lemma \ref{proj are aus-inj are bas}(2) and the dual argument to that of \cite[Proposition 4.9(1)]{BDGO21}. \cqfd

The following is the first main result of this section. It shows that $\B_C(R)$ (resp., $\A_C(R)$) is part of a right (left) AB-context. Its proof is similar to that of Theorem \ref{Gc-flat weak AB-cont} and Corollary \ref{Gc-flat approx}.

\begin{thm}Let $C$ be an $R$-module. The following assertions hold.
\begin{enumerate}
\item If $_RC$ is self-small and $\Sigma$-self-orthogonal, then the triple	$$\left(\Add_R(C),\cores(\widehat{\Add_R(C)}),\B_C(R)\right)$$ is a right weak AB-context in $R$-{\rm Mod}.  In this case,  for every  $N\in \cores (\widehat{\B_C(R)})$ there exist exact sequences $$0\to Y^N\to X^N\to N\to 0 \text{ and } 0\to N\to Y_N\to X_N\to 0$$ with $Y_N,Y^N\in \B_C(R)$ and $X_N,X^N\in \cores (\widehat{\Add_R(C)})$. 
\item  If $_RC$ is self-co-small and $\prod$-$\Tor$-orthogonal, then the triple $$\left(\A_C(R),\res(\widehat{\Prod_R(C^+)}),\Prod_R(C^+)\right)$$ is a left weak AB-context in {\rm Mod}-$R$.  In this case, for every $M\in \res(\widehat{\A_C(R)})$, there exist exact sequences of right $R$-modules $$0\to Y_M\to X_M\to M\to 0 \text{ and } 0\to M\to Y^M\to X^M\to 0$$ with  $X_M,X^M\in \A_C(R)$ and $Y_M,Y^M\in \res (\widehat{\Prod_R(C^+)})$.
\end{enumerate}
\end{thm}

\begin{rem}
Using all left and right AB-contexts constructed in this paper, one would use Lemmas \ref{left Frobnius} and \ref{right Frobnius}, together with \cite[Theorems B and C]{LY22}, to describe more left and right AB contexts in the category of complexes of $R$-modules, and in the category of representations of a left rooted quiver $Q$ with values in the category of $R$-modules.
\end{rem}

Now we investigate when the class $\Add_R(C)$ (resp., $\Prod_R(C^+)$) is the right (resp., left) hand class of a complete cotorsion pair.

\begin{prop}\label{special pc-proj preen}
Let $_RC$ be self-small, $\Sigma$-self-orthogonal and such that the canonical map $R\rightarrow\Hom_S(C,C)$ is an isomorphism. Then,  $\left( ^\perp\Add_R(C),\Add_R(C)\right) $ is a complete cotorsion pair if and only if $S$ is left perfect and right coherent, and $C_S$ is finitely presented and FP-injective.
\end{prop}
\proof First of all, notice that following the argument of \cite[Proposition 3.9]{KS98} (keeping in mind Lemma \ref{equiv of categ}(1)), one easily sees that the condition of $_RC$ being finitely generated can be substituted by that of $_RC$ being self-small.

$1.\Rightarrow 2.$ By hypothesis, the class $\Add_R(C)=(^\perp\Add_R(C))^\perp$  is closed under direct products and then $S$ is left perfect and right coherent, and $C_S$ is finitely presented  by \cite[Theorem 3.1 and Proposition 3.9]{KS98}. On the other hand,  we have $R^+\in (^\perp\Add_R(C))^\perp=\Add_R(C)$. Then, $R^+=C\otimes_S P$ for some projective $_SP$. Hence, 
$C^+\cong (R\otimes_R C)^+\cong \Hom_R(C,R^+)\cong \Hom_R(C,C\otimes_S P)\cong P$. Thus, $C_S$ is FP-injective by \cite[Theorem 1]{CS81}.

$2.\Rightarrow 1.$ First we prove that every $R$-module $M$ has a special $\P$-projective preenvelope.

By \cite[Theorem 3.1. and Proposition 3.9.]{KS98}, $\Add_R(C)$ is enveloping and by Lemma \ref{Wakamtsu Lemma} it suffices to show that each $\P$-envelope $f:M\to X$  is a monomorphism. If we prove that every injective module is $\P$-projective, then clearly the map $f$ will be a monomorphism. But since the class $\Add_R(C)$ is closed under direct summands, and using \cite[Theorem 3.1 and Proposition 3.9]{KS98} again we get that it is also closed under direct products, we only need to check that the injective cogenerator $R$-module $R^+$ is $\P$-projective.

By \cite[Theorem 1]{CS81}, $C^+$ is a flat $S$-module and then projective since $S$ is left perfect. So, by \cite[Theorem 3.2.11]{EJ00}, $R^+\cong \Hom_S(C,C)^+\cong C\otimes_S C^+$ is $\P$-projective. Now, it only remains to show that $(^\perp\Add_R(C))^\perp=\Add_R(C)$.

The inclusion $\Add_R(C)\subseteq (^\perp\Add_R(C))^\perp$  is clear. For the converse, let $X\in (^\perp\Add_R(C))^\perp$ and consider a special $\P$-projective preenvelope  of $X$: $0\to X\to F\to L\to 0$. Hence, $\Ext^1_R(L,X)=0$ and then the short exact sequence splits. Thus, $X$ is a direct summand of $F$ and so $X\in \Add_R(C)$.  \cqfd


To prove the dual of Proposition \ref{special pc-proj preen} we need the following lemma.

\begin{lem}\label{S is noetherian}
Assume that  $_RC$ is self-co-small and $C_S$ is finitely presented. Then, the following assertions are equivalent:
\begin{enumerate}
\item Every right $R$-module $M$ has an $\I$-injective (pre)cover.
\item  Every direct sum (limit) of $\I$-injective right $R$-modules is $\I$-injective.
\item $S$ is a right noetherian ring.
\end{enumerate}
\end{lem}
\proof $1.\Rightarrow 2.$  Consider a family $(X_i)_{i\in I}$ of $\I$-injective right $R$-modules  and let $M=\oplus_{i\in I} X_i$. By hypothesis, there exists an $\I$-injective precover $f:X\to M$. For each $i\in I$, there exists a morphism $g_i:X_i\to X$ such that $fg_i=\lambda_i:X_i\to M$ is the canonical injection. Now take the morphism $g:M\to X$ such that $g\lambda_i=g_i$ for each $i\in I$. Hence, $fg\lambda_i=fg_i=\lambda_i$ for each $i\in I$. Thus, $fg=1_{M}$ and then $f:X\to M$ is a split epimorphism. Hence, $M$ is a direct summand of $X$ which is $\I$-injective. Hence, $M=\oplus_{i\in I} X_i$ is $\I$-injective as well.

$3.\Rightarrow 2.$ Since $S$ is right noetherian, every direct limit of injective right $S$-modules is injective, so the implication follows from the canonical isomorphism $$\Hom_S(C,\varinjlim_{i\in I} E_i)\cong \varinjlim_{i\in I} \Hom_S(C,E_i)$$ for any direct family of injective right $S$-modules.

$2.\Rightarrow 1.$ For every set $I$, the right $R$-module $(C^+)^I\cong (C^{(I)})^+$ is pure-injective, so is every $\I$-injective. Then, $C^+$ is $\Sigma$-pure-injective by the hypotheses, that is, $(C^+)^{(K)}$ is pure-injective for every set $K$. Hence, $\Prod_R(C^+)$ is precovering by \cite[Proposition 6.10]{A2}. Moreover, by the hypotheses and \cite[Corollary 5.2.7]{EJ00} the class $\Prod_R(C^+)$ is covering.

$2.\Rightarrow 3.$ By \cite[Theorem 3.1.17]{EJ00} we need to prove that every direct sum of injective right $S$-modules is injective. But since every injective right $S$-module is a direct summand of some $(S^+)^I$, it suffices to show that any direct sum of copies of $(S^+)^I$ is an injective right $S$-module, for any set $I$.

Let $K$ be a set. Since $_RC$ is self-co-small, $((S^+)^I)^{(K)}=((\Hom_R(C,C)^+)^I)^{(K)}\cong ((C^+)^I\otimes_R C)^{(K)}\cong ((C^+)^I)^{(K)}\otimes_R C\in \inj(S)$ by Lemma \ref{equiv of categ}(2) since $(C^+)^I)^{(K)}$ is $\I$-injective by hypothesis. \cqfd

\begin{prop}
Assume that $_RC$ is self-co-small and $\prod$-$\Tor$-orthogonal, that $C_S$ is finitely presented and that the canonical map $R\rightarrow\Hom_S(C,C)$ is an isomorphism. Then, $\left( \Prod_R(C^+),\Prod_R(C^+)^\perp\right) $ is a  complete cotorsion pair if and only if $S$ is right noetherian and $C_S$ is injective.
\end{prop}
\proof Follow a dual argument to that of Proposition \ref{special pc-proj preen} using Lemma \ref{S is noetherian}. \cqfd

\begin{thm}\label{bass-mod-struc} (\textbf{The $\B_C$-Bass model structure})
Let $C$ be a semidualizing $(R,S)$-bimodule.  The following assertions are equivalent: 
\begin{enumerate}
\item  There exists a unique hereditary abelian model structure on $R$-\Mod such that:
		
$\bullet$ The fibrant objects coincide with the $\B_C$-Bass $R$-modules.
		
$\bullet$ The class of cofibrant objects coincides with  $^\perp \Add_R(C)$
		
$\bullet$ The trivially fibrant objects coincide with the $\P$-projective $R$-modules.
		
$\bullet$ The class of trivially cofibrant objects coincides with $^\perp \B_C(R)$.
		
\item $S$ is left perfect and right coherent, $C_S$ is FP-injective and $\Add_R(C)$ is closed under cokernels of monomorphisms. 
\end{enumerate}
	
In this case, the category $\A_{c,f}=\;^\perp\Add_R(C)\cap\B_C(R)$ is a Frobenius category. The projective-injective objects are exactly the $\P$-projective $R$-modules. Moreover, the homotopy category of the $\B_C$-Bass model structure is triangle equivalent to the stable category $\underline{^\perp\Add_R(C)\cap\B_C(R)}$. 
\end{thm}
\proof 
$1.\Rightarrow 2.$  It follows from the hypothesis that $(\Q,\tilde{\R})=\left(^\perp\Add_R(C),\Add_R(C)\right) $ is a  complete hereditary cotorsion pair. So, the assertion follows by Corollary \ref{special pc-proj preen}. 

The rest of this result can be proven as in Theorem \ref{GwF model structure} and  Corollary \ref{GcF-homotopy}. \cqfd

Dually, we have:

\begin{thm}\label{Aus-mod-struct} (\textbf{The $\A_C$-Auslander model structure})
Let $C$ be a semidualizing $(R,S)$-bimodule.  The following assertions are equivalent:
\begin{enumerate}		
\item  There exists a unique hereditary abelian model structure on {\rm Mod}-$R$ such that:
		
$\bullet$ The  cofibrant objects coincide with the   $A_C$-Auslander right $R$-modules.
		
$\bullet$ The class of fibrant objects coincide with $\Prod_R(C^+)^\perp$.
		
$\bullet$ The trivially cofibrant objects coincide with the $\I$-injective right $R$-modules.
		
$\bullet$ The trivially fibrant objects coincide with the $R$-modules in $\A_C(R)^\perp$.
		
\item  $S$ is right noetherian, $C_S$ is injective and  $\Prod_R(C^+)$ is closed under kernels of epimorphisms.
\end{enumerate} 	

In this case, the category $\A_{c,f}=\A_C(R)\cap \Prod_R(C^+)^\perp$ is a Frobenius category. The projective-injective objects exactly the $\I$-injective right $R$-modules. Moreover, the homotopy category of the $\A_C$-Auslander model structure is triangle equivalent to the stable category  $\underline{\A_C(R)\cap \Prod_R(C^+)}$.
\end{thm}

Let $C$ be a semidualizing $(R,R)$-bimodule. Recall \cite[Definition 5.1]{BEGO21a} that an extension closed full subcategory $\mathcal{F}$ of $R$-Mod is said to be a $C$-Frobenius category provided that it has enough projective and enough injective objects and $\mathcal{I}(\mathcal{F})=C\otimes_R\mathcal{P}(\mathcal{F})$ where $\mathcal{I}(\mathcal{F})$ and $\mathcal{P}(\mathcal{F})$  denote the full subcategories of the projective and injective objects of $\mathcal{F}$, respectively. In particular, $R$-Mod is $C$-Frobenius if and only if the injective $R$-modules coincide with the $\P$-projective modules. Note that $R$ is quasi-Frobenius if and only if $R$-Mod is $_RR$-Frobenius.

\begin{lem}
Assume that $C$ is a semidualizing $(R,R)$-bimodule. Then,  $R$-\Mod is $C$-Frobenius if and only if the projective $R$-modules coincide with the $\I$-injectives. 
\end{lem}
\proof  Straightforward. \cqfd




\begin{cor}
Assume that $C$ is a semidualizing $(R,R)$-bimodule. The following assertions are equivalent:
\begin{enumerate}
\item There exists a unique hereditary Hovey triple $\mathcal{M}_{\mathcal{I}}=(R{\rm -Mod},\,^\perp\B_C(R),\B_C(R)).$
\item $R$-\Mod is $C$-Frobenius. 
\item There is a unique hereditary Hovey triple $\mathcal{M}_{\mathcal{P}}=(\A_C(R),\A_C(R)^\perp,R{\rm -Mod}).$
\end{enumerate}
	
In this case,  the categories $\B_C(R)$ and $\A_C(R)$ are Frobenius categories whose projective-injective objects are exactly the $\P$-projective and $\I$-injective $R$-modules respectively. Moreover, we have the triangle equivalences $$\Ho(\mathcal{M}_{\mathcal{I}})\cong \underline{\B_C(R)} \text{ and } \Ho(\mathcal{M}_{\mathcal{P}})\cong \underline{\A_C(R)}.$$

\end{cor}

When $C=R$, the above two model structures coincide. Their homotopy category is known as the stable module category. The stable module category is the main object of study in modular representation theory and is discussed in detail by Hovey in \cite[Section 2.2]{Hov99}.

\begin{cor}[\cite{Hov99}, Section 2.2] 
The following assertions are equivalent:
\begin{enumerate}
\item There exists a unique hereditary Hovey triple $(R{\rm -Mod},\proj(R),R{\rm -Mod}).$
\item $R$ is quasi-Frobenius.
\item There exists a unique hereditary Hovey triple $(R{\rm -Mod} ,\inj(R),R{\rm -Mod})$.
\end{enumerate} 
\end{cor}

\bigskip
\noindent\textbf{Acknowledgment:} This work started when R. El Maaouy visited the mathematics department of the University of Almer\'ia as part of an Erasmus+ KA107 scholarship. He would like to thank the University of Almer\'ia for the warm hospitality and excellent working conditions. The authors J. R. Garc\'{\i}a Rozas and L. Oyonarte were partially supported by the grants by Ministerio de Econom\'{\i}a y Competitividad project of reference PID2020-113552GB-I00, and by the Junta de Andaluc\'{\i}a project with reference P20-00770.

\end{document}